\documentclass[12pt]{imsart}

\RequirePackage[OT1]{fontenc}
\RequirePackage{amsthm,amsmath,latexsym,amssymb}
\RequirePackage[round]{natbib}
\RequirePackage[colorlinks,citecolor=blue,urlcolor=blue]{hyperref}

\usepackage{graphicx}
 \usepackage{caption}
 \usepackage{subcaption}
\usepackage{amscd}
\usepackage{pdflscape}

\newcommand{\rasq}{R_\alpha^2}
\newcommand{\Ma}{M_\alpha}
\newcommand{\Mac}{M_{\alpha_c}}
\newcommand{\hatmat}{ P_n(\alpha)}
\newcommand{\hatmatAc}{ P_n(\alpha_c)}
\newcommand{\thtA}{{\boldsymbol \theta_\alpha}}
\newcommand{\scrptA}{\mathcal{A}}
\newcommand{\yn}{{\bf y}_n}
\newcommand{\Xn}{{\bf X}_n}
\newcommand{\mun}{{\boldsymbol \mu}_n}

\newcommand{\Xan}{{\bf X}_{\alpha}}

\newcommand{\Ba}{\boldsymbol{\beta}_\alpha}

\newcommand{\en}{{\bf e}_n}

\newcommand{\pna}{p(\alpha)}
\newcommand{\pnac}{p(\alpha_c)}
\newcommand{\sysq}{\mathcal{S}_y^{2}}
\newcommand{\sesq}{\mathcal{S}_e^{2}}

\arxiv{arXiv:0000.0000}

\startlocaldefs
\numberwithin{equation}{section}
\theoremstyle{plain}
\newtheorem{thm}{Theorem}[section]
\newtheorem{result}{Result}[section]
\newtheorem{remark}{Remark}[section]
\newtheorem{lemma}{Lemma}[section]
\endlocaldefs

\begin{document}

\begin{frontmatter}
\title{A Mixture of \lowercase{$g$}-priors for Variable Selection when the Number of Regressors Grows with the Sample Size }
\runtitle{Mixture of \lowercase{$g$}-priors for variable selection}
  
\begin{aug}
 \author{Minerva Mukhopadhyay\\
 Indian Statistical Institute \& Bethune College \\ \lowercase{{\small minervamukherjee@gmail.com}}}

\runauthor{Mukhopadhyay}



\end{aug}

\begin{abstract}
We consider variable selection problem in linear regression using mixture of $g$-priors. A number of mixtures are proposed in the literature which work well, especially when the number of regressors $p$ is fixed.
In this paper, we propose a mixture of $g$-priors suitable for the case when $p$ grows with the sample size $n$. We study the performance of the method based on the proposed prior when $p=O(n^b),~0<b<1$. Along with model selection consistency, we also investigate the performance of the proposed prior when the true model does not belong to the model space considered. We find conditions under which the proposed prior is consistent in appropriate sense when normal linear models are considered. Further, we consider the case with non-normal errors in the regression model and study the performance of the model selection procedure.
We also compare the performance of the proposed prior with that of several other mixtures available in the literature, both theoretically and using simulated data sets.
\end{abstract}



\end{frontmatter}

\section{Introduction}
We consider the regression problem with a response variable $y$ and a set of $p$ potential regressors $\mathnormal{x_1}, \mathnormal{x_2}, \ldots, \mathnormal{x_p}$. Let $\yn=\left(y_1, y_2, \ldotp\ldotp \right.$ $\left. \ldotp ,y_n \right)$ be a set of $n$ observations on $y$, and ${\bf X}_n=\left({\bf x}_1, {\bf x}_2, \ldots , {\bf x}_p \right) $ be the $n\times p $ design matrix, where ${\bf x}_i$ is the vector of $n$ observations on the $i^{th}$ regressor $\mathnormal{x_i}$, $i=1,2, \ldots, p$. We write
\begin{equation}
 \yn=\mun+\en, \label{1}
\end{equation}
where $\mun= E(\yn| {\bf X}_n)$ is the regression of $\yn$ on ${\bf X}_n$ and $\en$ is the vector of random errors.
If we assume the normal linear regression model, then $\mun=\beta_0 {\bf 1}+{\bf X}_n \boldsymbol{\beta}$ and $\en\sim N_n({\bf 0},\sigma^2 I)$. Here $\beta_0$ is the intercept, ${\bf 1}$ and ${\bf 0}$ are the $n\times 1$ vectors of ones and zeros, respectively, and $\boldsymbol{\beta}^{\prime}=\left(\beta_1,\beta_2, \ldots, \beta_p\right)$ is the vector of regression coefficients.

In this article, we study the variable selection problem. 
Given a set of $p$ available regressor variables, there are $2^p$ possible linear regression models depending on which regressors are included in the model. The space of all these models is denoted by $\scrptA$ and indexed by $\alpha$, where each $\alpha$ consists of a
subset of size $\pna~(0\leq\pna\leq p)$ of the set $\{1,2,\dots p\}$, indicating which regressors are selected in the model. The $\alpha^{th}$ model $\Ma$ is stated as
\begin{equation}
 \Ma~: \quad \mun=\beta_0 {\bf 1}+\Xan\Ba, \label{33} 
\end{equation}
 where $\Xan$ is a sub-matrix of ${\bf X}_n$ consisting of the $\pna$ columns specified by $\alpha$ and $\Ba$ is the corresponding vector of regression coefficients. We assume that all the components of $\Ba$ are non-zero. This ensures there is at most one true model in the model space. Our purpose is to choose the model $\alpha \in \scrptA$, which best explains the data. 

In a Bayesian approach, each model $\Ma$ is associated with a prior probability $p(\Ma)$ and the corresponding set of parameters ${\boldsymbol\theta}_\alpha=\left(\beta_0,\Ba,\sigma^2\right)$ involved in the model, is also associated with a prior distribution $p({\boldsymbol\theta}_\alpha| \Ma)$. Given the priors, one computes the posterior probability of $\Ma$ as
\begin{eqnarray}
   p(\Ma|\yn)= \displaystyle\frac{p(\Ma)m_\alpha(\yn)}{\sum_{\alpha \in \scrptA} p(\Ma)m_\alpha(\yn)},\label{40}
   \end{eqnarray}
where 
\begin{eqnarray}
\quad \quad m_\alpha(\yn) = \int p(\yn | \thtA, \Ma) p(\thtA|\Ma ) d\thtA 
\label{17}
\end{eqnarray}
is the marginal density of $\yn$ under $\Ma$ and $p(\yn | \thtA, \Ma)$ is the density of $\yn$ given $\thtA$ under $\Ma$.
In our search for a model, $p(\yn | \thtA, \Ma)$ will be taken to be normal. We consider the model selection procedure that selects the model in $\scrptA$ with the highest posterior probability.  

We consider a prevalent conventional prior on $\Ba$, the $g$-prior due to \cite{Zellner_1986}. Properties of $g$-prior are studied extensively in the literature (see, e.g., \cite{GF_2000}, \cite{Berger_Pericchi_2001}, \cite{bench_01}). The prior specification induced by the $g$-prior method crucially depends on the choice of the hyperparameter $g$ (see, e.g., \cite{Berger_Pericchi_2001}, \cite{Lng_et_al}). It has also been argued in the literature that this method is subject to inconsistencies like {\it Bartlett paradox} (see \cite{Bartlett_1957}, \cite{Jeffreys_1961}) and {\it information paradox} (see \cite{Zellner_1986}, \cite{Berger_Pericchi_2001}). \cite{Lng_et_al} considered a mixture on $g$ instead of considering a fixed $g$ to overcome these inconsistencies. Subsequently, a number of mixtures on $g$ are proposed in the literature. In this paper, we propose a mixture $\pi(g)$ on $g$,  
suitable for the case when $p$ grows with $n$.

We assume without loss of generality that the columns of $\Xan$ are centered, so that ${\bf 1}^{\prime}{\bf x}_i=0$ for all $i$.
 The intercept $\beta_0$ and the scale parameter $\sigma^2$ are common to all models and assumed to be independent of other parameters. We use standard non-informative priors for $\beta_0$ and $\sigma^2$. For justification of using such priors, see \citet[Sec 3.3]{Bayari_2012}. The vector $\Ba$ is assumed to follow a normal distribution with location parameter ${\bf 0}$ and scale $g \sigma^2 \left(\Xan^{\prime}\Xan \right)^{-1}$, where the hyperparameter $g$ has density $\pi(g)$. The complete prior specification is given by
\begin{eqnarray}
 &&p(\beta_0,\sigma^2 |\Ma)=\frac{1}{\sigma^2},~~
 \Ba|\left(\beta_0,\sigma^2,g,\Ma\right) \sim N_{p(\alpha)}({\bf 0},g\sigma^2 (\Xan^\prime\Xan)^{-1}) \label{4} \\
\mbox{and} && g\sim \pi(g).\nonumber  
\end{eqnarray}
We do not consider any specific prior probability on the model space. We only impose some conditions on model prior probabilities under which our results hold.
Similar setup has been considered by many authors, see, e.g., \cite{Lng_et_al}, \cite{Bayari_2012}.

Among the existing mixtures of $g$-priors, the earliest one, to the best of our knowledge, is due to \cite{Zellner_Siow}, who proposed a Cauchy prior on $\Ba$. Since Cauchy is an inverse gamma scale mixture of normal distributions, their prior proposition is considered as a mixture of $g$-priors. Other priors include hyper-$g$ and hyper-$g/n$ priors proposed by \cite{Lng_et_al}, generalized $g$-prior of \cite{George_2011} and robust prior proposed by \cite{Bayari_2012}. Henceforth, we will refer to these priors as Zellner-Siow prior, hyper-$g$ or hyper-$g/n$ prior, generalized $g$-prior and robust prior, respectively.

Inspite of existence of the above mixtures of $g$-priors there is not much discussion in the literature on which mixture one should use in a given situation. \cite{Bayari_2012} described some desirable properties a prior should satisfy in the context of Bayesian model selection. \cite{LS_2012} made an extensive simulation study to compare several priors. However, none of them considered the case when $p$ increases with $n$. \cite{George_2011} proposed a prior which is applicable when $p>n$, but proved consistency of their method for the case when $p$ is fixed. \cite{SC_2011} proved consistency for mixture of $g$-priors when $p$ increases with $n$ but their setup differs from the usual $g$-prior setup with respect to the covariance structure of the prior distribution of $\Ba$. \cite{WS_2014} investigated properties of different mixtures for the case with growing number of regressors. However, they only established results for Bayes factor consistency.

In this paper, we consider a scaled inverse chi-square prior on $g$ with appropriate parameters.
The Zellner-Siow prior belongs to this family. In a sense we consider a modified form of Zellner-Siow prior by choosing an appropriate scale parameter. An advantage of this prior is that it provides an approximation to the marginal density in (\ref{17}) with a closed-form expression which facilitates easy implementation and theoretical studies. Further, it satisfies many attractive consistency properties when $p$ increases with $n$. Most of the existing mixtures fail to be consistent for such rates of increase of $p$.
The good properties of this prior are not restricted to the case where the error distribution is normal. For a general class of error distributions with minimal assumptions, the proposed prior performs reasonably well in terms of consistency. 
In Section 2, we explicitly describe the form of the prior and the motivation for considering the same. 

In our investigation, we assume that the number of regressors $p$ increases with $n$ at a rate $p=O(n^b),~ 0<b<1$ and $p<n$. This is the so called {\it `large $p$ large $n$ regime'} and is of theoretical interest in contemporary research (see, e.g., \cite{FP_2004}, \cite{MGC_2010}, \cite{SKG_2012}, \cite{JR2012}). In practice, one can conveniently use methods applicable in {\it `large $p$ large $n$ regime'} when there is a sizable number of regressors compared to $n$ and the number of competing models in the model space is significantly large compared to $n$.   

We show that our proposed prior is consistent in appropriate sense for a large class of models under reasonable assumptions. We consider the following two cases separately. First, we consider the case when 
the true model belongs to the model space $\scrptA$, i.e., the true regression $\mun$ is as in (\ref{33}) for some $\alpha$. This case is referred to as the {\it `model true'} case. A well known notion of consistency in this regard is \emph{model selection consistency} which requires the posterior probability of the true model to go to one as $n\rightarrow\infty$. We examine model selection consistency of the proposed mixture of $g$, along with that of some other existing mixtures. 
We then consider the case when $\mun$ can be any unknown vector, not necessarily in the span of $\{{\bf 1}, {\bf x}_1, {\bf x}_2, \ldots , {\bf x}_p\}$. This case is referred to as the {\it `model false'} case. {\it `Model false'} case has been previously studied in \cite{Shao_1997}, \cite{CG_2006}, \cite{ARC_TS_2008} and  \cite{MSC_2013}. We investigate consistency of the proposed prior in this case using an appropriate notion of consistency.

The presence of information paradox in Zellner's $g$-prior remains one of the key motivations for considering mixture of $g$-priors. So, it is important to verify whether the proposed mixture can resolve the information paradox, i.e., is information consistent in the sense of \cite{Bayari_2012}. Along with the above notions of consistency we study information consistency as well. 

Sections 3 and 4 of the paper deal with model selection consistency.
In Section 3.1, we first consider the case when the error distribution is normal. 
In Section 3.2 we relax the condition of normality on $\en$. Here $e_i$s are assumed to be {\it i.i.d.} with mean 0 and finite fourth moment. 
In Section 4, we consider the performance of some other mixtures on $g$ with respect to model selection consistency in normal linear model setup. 
Section 5 deals with the {\it `model false'} case.
We find sufficient conditions under which the proposed prior is consistent in an appropriate sense under general error distributions.
In Section 6, we consider information consistency. 
In Section 7, we validate the performance of the proposed prior with extensive simulation studies. We study the performance of the proposed prior in comparison with several other existing priors in the literature. 
Finally, we make some concluding remarks in Section 8.
Proofs of all the main results are presented in the Appendix (Section 9).
\section{Scaled inverse chi-square mixture of $g$-priors}
We first motivate our proposal of a mixture on $g$. Most of the mixtures on $g$ in the existing literature 
are highly positively skewed having a unique modal point close to zero and a very flat decay. For example, hyper-$g$ and hyper-$g/n$ priors are $J$-shaped with modal point at $0$. Again, if we consider popular recommendations of $g$ in Zellner's $g$-prior, choices include the unit information prior ($g=n$, \cite{KR_1995}), the choice of $g$ related to the risk inflation criterion ($g=p^2$, see \cite{GF_1994}, \cite{GF_2000}), and the benchmark prior ($g=\max\{n,p^2\}$, \cite{bench_01}). Recently, \cite{MSC_2013} presented some theoretical results to explain why a relatively larger value of $g$ yields better results, especially when $p$ grows with $n$ and recommended using $g=n^2$ for practical purposes. From such recommendations, it seems reasonable to put relatively higher probability masses to higher values of $g$ for a mixture. Thus, there  persists a gap in the domain of $g$ getting relatively higher mass when a fixed $g$ is 
considered compared to that of a mixture. 
We now propose a mixture which gives more probability mass to a range of relatively higher values of $g$ compared to the existing mixtures.
We consider the {\it scaled inverse chi-square} mixture $\pi(g)$ on $g$ with scale parameter $\tau^2=n^2$ and degrees of freedom $\nu$, given by
\begin{equation}
\pi(g)=\frac{\left( \tau^2 \nu/2\right)^{\nu/2}}{\Gamma{(\nu/2)}}\frac{\exp\left[- \tau^2\nu/(2g) \right]}{g^{1+\nu/2}}, \hspace{.1 in} g>0,~\nu>0, ~\tau^2>0.\label{5} 
\end{equation}
Although, the hyperparameter $\nu$ can take any positive value, we recommend using values between $1$ and $p$. In this paper, we will consider two extreme choices of $\nu$, namely, $\nu=1$ and $\nu=p$. Note that such choices of hyperparameters ensure that the prior has a unique mode at $ n^2 \nu /(\nu+2)$ and a very flat decay.

The use of inverse gamma distribution in a scale mixture of normal priors for $\Ba$ is a common practice in Bayesian model selection. It has already been stated that the Zellner-Siow prior for $\Ba$ is an inverse gamma scale mixture of normals with shape parameter $1/2$ and scale parameter $n/2$. In the context of linear regression models with shrinkage priors, \cite{PC_2008} and \cite{Hans_2009}, while introducing Bayesian version of lasso, 
used inverse-gamma priors for similar normal scale mixtures for $\Ba$. 
The proposed prior in (\ref{5}) is same as the inverse gamma prior with shape parameter $\nu/2$ and scale parameter $n^2\nu/2$. 

An advantage of considering this prior is that it yields a closed form approximation to the marginal density, which is similar to the form of the marginal for a $g$-prior with some fixed choice of $g$. Availability of closed form marginals (posterior probabilities) is not necessary for good inference but it serves as a desirable property for easy implementation. 
This approximation makes the application of this prior simple and theoretically tractable.
\subsection{Posterior probability}
For the linear model setup, the vector of parameters in model $\Ma$ is given by ${\boldsymbol \theta}_{\alpha}=\left(\beta_0,\Ba,\sigma^2,g\right)$. The marginal density $m_\alpha(\yn)$ in (\ref{17}) is given by 
\begin{eqnarray*}
 && m_\alpha(\yn)\\
 &&=\int p(\yn |  \beta_0, \Ba, \sigma^2, \Ma) p(\Ba|\beta_0, \sigma^2, g,\Ma ) p(\beta_0 , \sigma^2) \pi(g) d(\beta_0,\Ba,\sigma^2,g), 
\end{eqnarray*}
where $p(\yn | \beta_0, \Ba, \sigma^2, \Ma)$ is p.d.f. of the $n$-variate normal distribution with mean $\beta_0 {\bf 1} + \Xan\Ba$, dispersion matrix $\sigma^2 I$ and $p(\Ba|\beta_0, \sigma^2, g,\Ma )$, $p(\beta_0 , \sigma^2)$, $\pi(g)$ are as in equations (\ref{4}) and (\ref{5}).

Integrating the integrand above with respect to $\beta_0$, $\Ba$ and $\sigma^2$, we get a closed form expression which leads to
\begin{align}
 \label{8} m_\alpha(\yn)=\frac{\Gamma {(n-1)/2}}{\pi^{(n-1)/2}\sqrt{n}} \left(\sysq \right)^{-(n-1)/2} \int_{0}^{\infty} \frac{(1+g)^{(n-1-\pna)/2}}{\left[1+g(1-\rasq) \right]^{(n-1)/2}}~ \pi(g) dg,
\end{align}
where $\sysq=\|\yn-\bar{y}_n {\bf 1} \|^2/n$, $(1-\rasq)= \yn^{\prime} \left(I-\hatmat\right)\yn/(n\sysq)$, $\hatmat=Z_{n\alpha}\left(Z_{n\alpha}^{\prime} Z_{n\alpha} \right)^{-1}Z_{n\alpha}^{\prime}$, 
and $Z_{n\alpha}=\left({\bf 1}~~\Xan\right)$.

Note that the marginal density of the intercept only model $M_{N}~:~\yn=\beta_0 {\bf 1}+\en$, which will be referred to as the {\it null} model, does not involve the hyperparameter $g$. It can be obtained as a special case of the marginal in expression (\ref{8}) by putting $\rasq=0$ and $\pna=0$. 

For models $\alpha\in\scrptA\smallsetminus \{N\}$, the marginal density given by the proposed prior (\ref{5}) does not have a closed form. But, we can make an approximation of the marginal density in a closed form expression when the proposed mixture (\ref{5}) is used in (\ref{8}).  
When $p$ is fixed, this approximation can obtain an accuracy of order $n^{-1}$ with probability tending to 1, which is the same as the accuracy of the Laplace approximation for fixed $p$ (see \cite{KR_1995}). When $p$ is increasing, the Laplace approximation may not be valid for the integral in (\ref{8}) for commonly used priors on $g$, since the integrand may not be {\emph Laplace~regular} (see \cite{KTK_1990}).
When $p=O(n^b),~0<b<1$ and $\nu=1$ (or, when $\nu$ is free of $n$), the approximation is less accurate with an error of the order $n^{-(1-b)}$. But, if $\nu=p$ (or, if $\nu$ is of same order of $n$ as $p$), the approximation still attains an accuracy of the order $n^{-1}$.

We first state the assumptions under which the approximation holds. 
Throughout this paper $\yn$ is modeled as (\ref{1}) and we assume the following:

\vskip5pt
(A.1)~ $\mun^\prime \mun/n < M$ for some constant $M>0$ as  $n\rightarrow \infty$.  
\vskip5pt

We make a mild assumption on distribution of $\en$.

\vskip5pt
(A.2)~The errors $e_1, e_2, \ldots, e_n$ are {\it i.i.d.} with a common density having mean $0$ and finite fourth order moment. 
\vskip5pt  

We now state the result.
\begin{result}\label{res_1}
  Consider the set of priors (\ref{4}) and (\ref{5}) with $\nu$ varying from 1 to $p$. Under assumptions (A.1) and (A.2), the marginal density in (\ref{8}) satisfies the following:
\begin{align*}
 && \quad \quad \quad m_{\alpha}(\yn) \leq \widetilde{m}_{\alpha}(\yn) \left(1+\frac{p}{\nu n}O(1)\right) \quad \quad \quad \quad \quad \quad \quad \quad \quad \quad \quad \quad \quad \quad \\
\mbox{and} && \quad \quad \quad m_{\alpha}(\yn) \geq \widetilde{m}_{\alpha}(\yn) \left(1+\frac{p}{\nu n}Op(1)\right),\quad \quad \quad \quad \quad \quad \quad \quad \quad \quad \quad \quad \quad
\end{align*}
where 
\begin{align}
&&\widetilde{m}_{\alpha}(\yn)=\hspace{4.3 in}\nonumber \\
&&~\frac{\Gamma {((n-1)/2)}\Gamma((\nu+\pna)/2)}{\sqrt{n}\Gamma(\nu/2)} \left(\pi \sysq \left(1-\rasq\right) \right)^{-(n-1)/2} \left(\frac{n^2\nu}{2} \right)^{-\pna/2},\nonumber
\end{align}
uniformly in $\alpha$, for any $\alpha\in \scrptA\smallsetminus \{N\}$ as $n\rightarrow\infty$.
\end{result}
\vskip5pt

From Result \ref{res_1}, we find an approximation to the marginal density as
\begin{equation}
  m_{\alpha}(\yn) \approx \widetilde{m}_{\alpha}(\yn), \label{18}
\end{equation}
in the sense that the ratio of $m_{\alpha}(\yn)$ and $\widetilde{m}_{\alpha}(\yn)$ goes to 1 in probability. This approximation holds uniformly in $\alpha$ since the $O(1)$ and $O_p(1)$ terms can be made free of $\alpha$. 
It is easy to check from Result \ref{res_1} that if $\nu=n^r$ for some $0\leq r\leq b$, then the approximation is accurate upto an order $1/n^{1+r-b}$. A simulation study on the performance of this approximation is added in the supplementary file.

\section{Model selection consistency of scaled inverse chi-square mixture}
In this section, we assume that the true mean $\mun$ can be expressed as a linear combination of a subset of the $p$ regressors. Let $\Mac$, ${\alpha}_c\in\scrptA$ be the true model. An ideal model selection procedure should identify the true model in this framework. Therefore, a natural criterion in this case is model selection consistency,  
which is achieved if posterior probability of the true model, given by (\ref{40}), converges to $one$ in probability, i.e.,
\begin{equation}
  p(M_{\alpha_c}|\yn) \xrightarrow{p} 1 \quad\mbox{as} \quad n\rightarrow \infty. \label{12}
\end{equation}

In Section 3.1, we provide sufficient conditions under which (\ref{12}) holds when the error distribution in (\ref{1}) is normal. In Section 3.2, we relax the assumption of normality and assume that $\en$ follows any distribution satisfying assumption (A.2). In this wider class of distributions, we investigate how the set of sufficient conditions for achieving (\ref{12}) modify under the nested model setup. 

\subsection{Model selection consistency with normality assumption}
Throughout this subsection, we will assume that  $\en\sim N_n\left( {\bf 0},\sigma^2 I \right)$. We show that when $p$ grows with $n$ then under appropriate conditions, model selection consistency is achieved by the proposed prior considering all $2^p$ models in model space. 

We split the model space into three mutually exclusive and exhaustive parts as follows: 

\noindent$\scrptA_1=\{\alpha\in\scrptA~:~\Ma\supset\Mac,$ $\alpha\neq\alpha_c\}, \;\scrptA_2=\{\alpha\in\scrptA~:~\alpha \notin \scrptA_1, \alpha\neq\alpha_c\}$ and $\;\{\alpha_c\}$ where $M_{\alpha_c}$ is the true model. 
We assume that

\vskip5pt
 (A.3)~  $\varliminf_{n\rightarrow\infty} n^s \min_{\alpha \in \scrptA_2}~ \mun^\prime(I-\hatmat) \mun/n > \delta$ for some constants~~$\delta>0$ ~~and~~ $0\leq s < 1$.\\
We impose a general restiction on model prior probability as
 
 \vskip5pt
(A.4)~$ \displaystyle\max_{\alpha, \alpha^\prime \in \scrptA} p(\Ma)/p(M_{\alpha^\prime}) \leq C$  \hspace{.1 in} for some constant~ $C>0$.
\vskip5pt

A remark on each of the assumptions is made below (see Remark \ref{rm_7} and Remark \ref{rm_8}).
\begin{thm}
\label{thm_1}
 Let $\yn$ be as in (1.1) with $\mun$ satisfying (A.1) and $\en\sim N_n\left({\bf 0},\sigma^2 I\right)$.  If $p=O(n^b)$, then the prior specification given by (\ref{4}) and (\ref{5}) is model selection consistent for $0<b<2/5$ when $\nu=1$, and $0<b<1/2$ when $\nu=p$, provided (A.3) holds with $s<(1-b)/2$ and (A.4) holds.
\end{thm}
 \begin{remark} \label{rm_6}
  This result is different from those obtained by \cite{WS_2014} who have shown Bayes factor consistency or pairwise consistency for growing number of regressors. Our result deals with the asymptotic behavior of the posterior probability of the true model considering all the $2^p$ models in the model space, and it is indeed a much stronger result than pairwise consistency. See in this context \citet[p.652]{JR2012}.
 \end{remark}
 \begin{remark}\label{rm_7}
  Assumption (A.3) for $s=0$ was assumed by many authors (see, e.g., \cite{bench_01}, \cite{Lng_et_al}, \cite{Bayari_2012}). It is a key assumption for model selection consistency which ensures that the models can be differentiated. Here we relax this assumption by allowing $s>0$, which is a natural extension for the situation when $p$ grows with $n$. 
 \end{remark}
\begin{remark} \label{rm_8}
Assumption (A.4) may seem a bit restrictive when $p$ grows with $n$. For the results in this paper to hold we do not actually need assumption (A.4); we can use much weaker versions of the assumption. For each of these results we mention the weaker version of the assumption needed for the proof. As a whole we work with assumption (A.4) for simplicity in presentation.

For Theorem \ref{thm_1} to hold, we only need $ \max_{\alpha \in \scrptA} p(\Ma)/p(\Mac) \leq C$ for some constant $C>0$. This is a reasonable assumption, since this only indicates that the true model may not have a prior probability arbitrarily close to zero, which is necessary to achieve consistency.
\end{remark}

\subsection{Model selection consistency in general settings} \label{S1}
In this subsection, we extend our results to situations where the error distributions belong to a larger class satisfying assumption (A.2).
We investigate the strength of the model selection algorithm when the distribution of the errors is non-normal and the same model selection rule (based on the normal likelihood) is used. In other words we investigate robustness of our model selection rule for non-normal errors.
Unlike the case for normal errors, here we do not consider all the $2^p$ models. Rather, we restrict our search within a class of nested models. By nested models, we mean a set of models where for every pair of models, say, $M_1$ and $M_2$, either $M_2\subset M_1$, or $M_1\subset M_2$.
So, the index set $\scrptA^{*}$ of all nested models can be expressed as 
\begin{equation*}
\scrptA^{*}=\left\{\{\phi\},\{1\},\{1,2\},\ldots,\{1,2,\ldots,p\} \right\}, \quad \mbox{with} \quad\scrptA^{*} \subset \scrptA.
\end{equation*}
Note that $\scrptA^*$ has $p+1$ different models. When $p=O(n^b)$ with $0<b<1$, the number of models in $\scrptA^{*}$ also increases with $n$. While the cardinality of $\scrptA$ is exponential in $p$ (i.e., $2^p$), for $\scrptA^{*}$ it is linear in $p$.
On one hand generalization to non-normal errors broadens the scope of the model selection algorithm, on the other, restriction to the class of nested models abridges the model space.

Comparison of the set of nested models is of great practical interest. The situation with a model space like $\scrptA^{*}$ may occur, for example, when we have information on relative importance of the regressors and the regressors can be ordered accordingly. 
Model selection in nested models has been widely studied in the Bayesian paradigm when the error distribution is normal (see, e.g., \cite{D1992}, \cite{M1997}, \cite{CG_2008}, \cite{WS_2014}). Unlike these authors, who study Bayes factor consistency, we consider model selection consistency restricted to the space of nested models when $p$ is of the said order. We summarize our findings in the following theorem.
\begin{thm}
 \label{thm_3}
  {\it Let $\yn$ be as in (1.1) with $\mun$ satisfying (A.1) and $\en$ satisfying (A.2).  If $p=O(n^b)$, then the prior specification given by (\ref{4}) and (\ref{5}) with $\nu=1$ or $\nu=p$ is model selection consistent in $\scrptA^*$ for any $0<b<1$ provided (A.3) holds with $s<(1-b)/2$ and (A.4) holds.}
\end{thm}
\begin{remark}\label{rm_9}
 Here also we do not need assumption (A.4) to hold strictly. For Theorem \ref{thm_3} 
to hold we only need $ \max_{\alpha \in \scrptA^*} p(\Ma)/p(\Mac) \leq C\sqrt{n}$ for some constant $C>0$. This assumption is quite general and includes many of the popular class of model prior probabilities.  
\end{remark}

 \section{Properties of Some Existing Mixtures}
 We now investigate the performance of some other mixtures from the perspective of model selection consistency. The beta prime (beta of second kind) prior is the commonly used prior for $g$ (see \cite{Lng_et_al}, \cite{George_2011}, \cite{Bayari_2012}). Therefore, it is worth investigating the performance of this prior when $p$ grows with $n$. Let $g$ follow a beta prime distribution with parameters $\gamma_0$ and $\gamma_1$, then
\begin{equation}
\pi(g)=\frac{\Gamma(\gamma_0+\gamma_1)}{\Gamma(\gamma_0) \Gamma(\gamma_1)} g^{\gamma_0-1}(1+g)^{-(\gamma_0+\gamma_1)}, \quad 0< g<\infty,~ \gamma_0>0,~\gamma_1>0. \label{14}
\end{equation}
We show that when $p=n^b$, $0<b<1$, then for some inappropriate specification of the hyperparameters $\gamma_0$ and $\gamma_1$, the model selection rule given by the set of priors (\ref{4}) and (\ref{14}) becomes inconsistent. 

\begin{thm}
 \label{thm_2}
Consider the setup of Theorem \ref{thm_1} and let the true model be the null model. If the number of regressors $p=n^b$, $0<b<1$, then the set of priors given by (\ref{4}) and (\ref{14}) with $\gamma_1>\epsilon$ for some $\epsilon >0$ free of $n$, is inconsistent if $(\gamma_0/\gamma_1) = O(n^{2b})$.
\end{thm}

\begin{remark}
 \label{rm_1} 
In hyper-$g$ prior $\gamma_0=1$ and $\gamma_1=(a/2-1)$ for some $a>2$ free of $n$. Hence, it follows from the above theorem that hyper-$g$ prior is inconsistent for any $b>0$. It is already shown in \cite{Lng_et_al} that hyper-$g$ prior is not consistent under the null model even for fixed $p$. Hyper-$g/n$ prior remains consistent when $p$ is fixed, but it fails to be consistent if $p=n^b$, for any $b>0$ (the proof is in the supplementary file). 
\end{remark}
\begin{remark}
 \label{rm_2}
 Generalized $g$-prior has $\gamma_0=A+1$ and $\gamma_1=\mathnormal{B}+1$ where the authors recommend using $A=(n-\pna-1)/2-B$ and $B<1/2$ for the case when $p<n$. Hence, it is easy to check that generalized $g$-prior is inconsistent for this recommended settings if $b\geq1/2$. 
 \end{remark}
 \begin{remark}
 \label{rm_3} 
 The robust prior can also be expressed as a truncated scaled beta prime distribution as $(g+B)/(\rho_{\alpha}(n+B))-1\sim beta~prime(1,A)$ where $A>0,~B>0,~\rho_{\alpha}>B/(B+n)$. The recommended choices of hyperparameters are $A=1/2,~B=1$ and $\rho_{\alpha}=1/(1+\pna)$. It has also been recommended that $\rho_{\alpha}$ should be free of $n$. This makes choice of the parameter $\rho_\alpha$ difficult when $p=n^b$, since in that case the choice of $\rho_{\alpha}$ involves $n$. We check with two choices of $\rho_\alpha$, a constant $\rho_\alpha$ and $\rho_\alpha=1/(1+\pna)$. It has been shown in the supplementary file that when $p=n^b$ a necessary condition for consistency of the robust prior under the null model is $b<1/2$, for both choices of $\rho_\alpha$.
 
 There are some other beta shrinkage priors available in the literature. A list of which can be obtained in \cite{LS_2012}. Similar results can be derived for them.
 \end{remark}
 \begin{remark}\label{rm_4}
  As we have already mentioned, Zellner Siow prior is also an inverse gamma prior with scale $n$, whereas the proposed one has scale $n^2 \nu/2$. It can be shown that a sufficient condition of Zellner-Siow prior to be consistent under null model is $b<1/5$ (the proof is similar to {\it Case II} in proof of Theorem \ref{thm_1} in Section \ref{Sect1}) . The increment of scale from the order of $n$ to $n^2$ makes the prior more reasonable for {\it `large $p$ large $n$ regime'} which is reflected in the increment of $b$ from $1/5$ to $1/2$. 
  
  However, similar improvement is not expected from all the priors we mentioned before. For example, if we change the scale of hyper-$g/n$ prior from $n$ to $n^r$ for any $r>1$ free of $n$, it still remain inconsistent when the null model is true, and the number of regressors is $p=n^b$, for any $0<b<1$. The proof is similar in idea to the proof of the result stated in Remark \ref{rm_1}.
   \end{remark}
   \begin{remark}
    \label{rm_5}
    As mentioned in Remark \ref{rm_8}, here also we do not need a strong assumption like (A.4). For all the above results in this section to hold, we only need to assume
   
   (A.4*)~~for all $ \alpha \in \scrptA,\quad p(\Ma)/p(M_{N}) \geq \delta^{\pna}~$ for some constant~ $\delta>0$. 
   
In situations when $M_{N}$ is a candidate for the true model, one may like to put additional penalty to more complex models. Therefore, the prior probabilities for the high-dimensional models may be quite small compared to that of low dimensional models. Assumption (A.4*) gives us the scope to consider those set of prior probabilities also. Consider, for example, the case where the prior inclusion probability of each regressor is $q$, which leads to the Bernoulli prior probability $q^{\pna} (1-q)^{p-\pna}$ to the model $\Ma$. The Bernoulli prior has been used by many authors (see, e.g.,\cite{GMc_1993}). It is easy to see that this prior setting satisfies assumption (A.4*).   
   \end{remark}
  Note that the conditions on $b$ mentioned in Theorem \ref{thm_1} are \emph{sufficient} to achieve consistency for the proposed prior. The condition in Remark \ref{rm_4} is also sufficient for Zellner-Siow prior, whereas the above conditions on beta shrinkage priors (i.e., hyper-$g/n$ prior, generalized $g$-prior and robust prior) are only \emph{necessary} to achieve consistency. The range of $b$ sufficient to achieve consistency is essentially a subset of the range of $b$ necessary to achieve consistency. Also, achieving consistency under the null model is only a part of achieving full posterior consistency. Sufficient conditions for full posterior consistency may even be stronger.


\section{Consistency of Scaled Inverse Chi-square Prior in `Model False' Case }
In Sections 3 and 4, we have considered situations when the true mean $\mun$ in (\ref{1}) belongs to the span of $\{ {\bf 1}, {\bf x}_1, \ldots, {\bf x}_n\}$. We will now consider a more general scenario where $\mun$ is any $n$-dimensional vector, i.e., the true model does not necessarily belong to the model space $\scrptA$. 
Several authors have studied related problems of linear model selection under this framework (see, e.g., \cite{Li_1987}, \cite{Shao_1997}, \cite{CG_2006}, \cite{ARC_TS_2008}, \cite{MSC_2013}). 
Of course, one cannot compute the posterior probability of the true model here, and the usual notion of model selection consistency cannot be used in this scenario. To validate a model selection rule we therefore adopt an alternative notion of consistency suited for the {\it model false} case as used in \cite{MSC_2013}. 

Here consistency of a model selection procedure refers to the property of choosing the model which is closest to the unknown true model among all candidate models in the model space $\scrptA$ (in an asymptotic sense).
We consider the Kullback-Leibler divergence as the measure of distance between two probability distributions. 
We define the distance $\Delta_n(\alpha)$ between the true distribution (of $\yn$) and the model $\Ma$ as the minimum of the Kullback-Leibler distance 
with respect to the underlying parameters $(\beta_0,\Ba)$ of the distribution under $\Ma$. 
One would naturally like to choose a model $M_{\alpha^*}$ which is as close as possible to the true distribution (i.e., for which $\Delta_n(\alpha^*)=\min_{\alpha\in\scrptA} \Delta_n(\alpha)$). One can find the model $M_{\alpha^*}$ only if the true distribution were known, which is not the case here. 
We show that our model selection procedure chooses a model which is closest to the unknown true model in an asymptotic sense.

We consider a general class of error distributions satisfying assumption (A.2) and the whole set of $2^p$ models $\scrptA$ for comparison.
We make the following assumption which is analogous to assumption (A.3), by replacing $\scrptA_2$ in (A.3) by $\scrptA$ as

\vskip5pt
(A.3*) ~$\varliminf_{n\rightarrow\infty} n^s \min_{\alpha \in \scrptA} \mun^\prime(I-\hatmat) \mun /n > \delta \hspace{.1 in}$ for some constant ~~ $\delta>0$ ~~and~~ $0\leq s < 1$.
\vskip5pt 
\noindent Note that in the {\it model false} case $\scrptA_2=\scrptA$.
 
 Let the true distribution of $\yn$ has a density function $f$. It can be easily verified that the Kullback-Leibler distance between the true distribution given by the density function $f$, and the distribution $ N\left({\bf 1} \beta_0+\Xan\Ba, \sigma^2 I \right)$ under $\Ma$ equals  
\begin{eqnarray*}
 &&\int f\left( \yn \right) \log f\left( \yn \right) d\yn +\frac{n}{2}\left(1+\log\sigma^2 \right)\\
 && \quad \quad \quad +\frac{1}{\sigma^2} \left(\mun-{\bf 1} \beta_0-\Xan\Ba\right)^\prime \left(\mun-{\bf 1} \beta_0-\Xan\Ba\right).
\end{eqnarray*}
The distance $\Delta_n(\alpha)$ between the true model $f$, and the model $\Ma$ is obtained by minimizing the above with respect to $(\beta_0, \Ba)$, as follows
\begin{eqnarray}
 \Delta_n(\alpha)
 &=&\int f\left( \yn \right) \log f\left( \yn \right) d\yn +\frac{n}{2}\left(1+\log\sigma^2 \right)+ D_n(\alpha), \label{23}
\end{eqnarray}
where
 \begin{equation}
D_n(\alpha)=\frac{1}{2\sigma^2} \mun^{\prime}\left(I-\hatmat\right)\mun. \label{20}
\end{equation}
Note that the first two terms of $\Delta_n(\alpha) $ in (\ref{23}) do not involve $\alpha$. Therefore, $\mbox{ argmin}_{\alpha} \Delta_n(\alpha)=\mbox{ argmin}_{\alpha} D_n(\alpha)$ and we show that
\begin{equation}
 \frac{D_n(\hat{\alpha})}{\min_{\alpha\in\scrptA}D_n(\alpha)} \xrightarrow{p}1 \mbox{~~as~~}n\rightarrow\infty,
 \label{25}
\end{equation}  
where $M_{\hat{\alpha}}$ is the model chosen by the model selection rule based on the proposed prior. 

\begin{thm}
 \label{thm_5}
{\it Let $\yn$ be as in (\ref{1}) with $\mun$ being any real vector in $\mathbb{R}^n$ satisfying (A.1) and $\en$ satisfying (A.2). Suppose the assumption (A.4) holds and (A.3*) holds with $s<(1-b)/2$. If the number of regressors $p=O(n^b)$, $0<b<1$, then the set of prior (\ref{4}) and (\ref{5}) is consistent in the sense that (\ref{25}) holds for any $0<b<1$ when $\nu=1$ or $\nu=p$.} 
\end{thm}

In Theorem \ref{thm_5}, 
it is only assumed that $\yn$ is the sum of two components, namely, the true mean $\mun$ and the random error $\en$. Here, $\mun$ is allowed to be arbitrary and $\en$ can follow any distribution satisfying assumption (A.2); even symmetry or continuity is not required. Thus, given the additive structure of $\yn$ as in (\ref{1}), consistency is obtained in a much general setting.

Lastly, assumption (A.4) can be relaxed to a great extent here. For Theorem \ref{thm_5} to hold, we only need the following

\vskip5pt
(A.4*) $\max_{\alpha,\alpha^\prime\in\scrptA} p(\Ma)/p(M_{\alpha^\prime}) \leq C n^r$ for some $C>0$ and $r>0$.

\vskip5pt
\noindent This assumption is satisfied for a very large class of prior probabilities on the model space.
 \section{Information Consistency}
   The criterion of information consistency is considered by several authors (see, e.g., \cite{Jeffreys_1961}, \cite{Berger_Pericchi_2001}, \cite{BG_2008}, \cite{Lng_et_al}, \cite{Bayari_2012}). 
 While comparing the null model with any model $\Ma$, suppose that $\| \hat{{\boldsymbol \beta}}_\alpha \|^2 \rightarrow\infty$ (or equivalently, the usual $F$ statistics goes to $\infty$) with both $n$ and $\pna$ are fixed, $\hat{{\boldsymbol \beta}}_\alpha $ being the least squares estimator of $\Ba$. This is considered as a very strong evidence supporting the model $\Ma$, and it is expected that the Bayes factor for comparing model $\Ma$ to the null model would go to $\infty$. The property that the Bayes factor goes to $\infty$ whenever $\| \hat{{\boldsymbol \beta}}_\alpha \|^2 \rightarrow\infty$ with fixed $n$ and $\pna$ is termed as information consistency in \cite{Bayari_2012}. However, this does not hold in the case of Zellner's $g$-prior. For mixture of $g$-priors, \citet[Theorem 2]{Lng_et_al} give a sufficient condition which ensures information consistency.
The following result gives conditions under which the proposed mixture in (\ref{5}) is information consistent. 
\begin{result}
 \label{res_2}
 Consider the set of prior probabilities (\ref{4}). Then the mixture on $g$ given by (\ref{5}) is information consistent if $n\geq p+1$ when $\nu=1$ and if $n\geq 2p$ when $\nu=p$. 
\end{result}
The proof of this result is in the supplementary file.

 Note that for $\nu=1$ the proposed prior is information consistent with minimal sample size, i.e., information consistency holds for any $n>p$ (see \cite{Lng_et_al} in this context). But for $\nu=p$, the proposed prior fails to be information consistent with minimal sample size. 

 \section{Performance of The Proposed Prior on Simulated Datasets}
  In this section we validate the performance of the proposed prior using simulated datasets. We present simulation results for model selection consistency under different simulation schemes. In each case, we consider our proposed prior with two choices of the hyperparameter $\nu$, {\it viz.}, $\nu=1$ ({\emph proposed prior I}) and $\nu=p$ ({\emph proposed prior II}). Along with the proposed prior we also consider four other priors on $g$, namely, Zellner-Siow prior, hyper-$g/n$ prior, generalized $g$-prior and robust prior.
  
Our results are designed for the case when $p$ increases with $n$, and therefore, we consider moderately large $p$ compared to $n$. Three choices of $n$ ($n=50,100,150$) and two choices of $p$ ($p+1=30,50$) for each $n$ have been considered.
 
 The theoretical results are not confined to the case with normal errors; any error distribution satisfying assumption (A.2) can be considered. We consider three different error distributions, namely, normal, Laplace and $t$ with degrees of freedom 3 ($t_{(3)}$). Note that $t_{(3)}$ does not satisfy the fourth order moment condition of (A.2). Moreover, we consider all the $2^p$ models in case of Laplace and $t_{(3)}$ distributions also, although our theoretical results for general case allow only nested model setup. We consider these settings to check the performance of the proposed mixture when some of the assumptions of Theorem \ref{thm_3} do not hold. The simulation scheme is described as follows.
 
 For each combination of $(n, p)$, we generate $n$ values of each of the $p$ regressors $x_1, x_2, \ldots , x_p$ and this gives the full design matrix $\Xn$. We choose $p$ numbers $\xi_i, i=1, \ldots , p$ and generate the $n$ values of the $i^{th}$ regressor $x_i$ from an $N(\xi_i, 1)$ distribution, $i=1, \ldots , p$. We assume that the $n$ values of the $i^{th}$ regressor are coming from a homogeneous population. In order to fix a ``true" model, we choose its dimension $p(\alpha_c)$ and then choose the $p(\alpha_c)$ non-zero regression coefficients $\beta_i$'s, the intercept $\beta_0$ in the true model and also a value for the error variance $\sigma^2$. The $p(\alpha_c)$ columns of the design matrix ${\bf X}_{\alpha_c}$ for the true model are chosen at random from the $p$ columns of $\Xn$.
Here, $(\xi_1, \ldots , \xi_p)$ is chosen as a random permutation of $(0.2, 0.4, \ldots , 0.2 \times p )$. The dimension of the true model $p(\alpha_c)$ is chosen as $[p/2]$ and the $p(\alpha_c)$ non-zero regression coefficients $\beta_j$'s and the intercept $\beta_0$ in the true model are randomly chosen from the set $\{ -0.2, 0.4, \ldots , (-1)^p \times 0.2 \times p \}$. Lastly, we choose $\sigma = 1$. 

After choosing the dimension $p(\alpha_c)$, the coefficients $\left( \beta_0,{\boldsymbol\beta}_{\alpha_c}\right)$, the error variance $\sigma^2$ of the true model and the design matrix $\Xn$, we generate $\en$ from normal, Laplace and $t_{(3)}$ distributions, each with location vector ${\bf 0}$ and dispersion matrix $\sigma^2 I$. The vector of observations $\yn$ is obtained by adding $\mun={\bf 1}\beta_0+ {\bf X}_{n\alpha_c} \boldsymbol{\beta}_{\alpha_c}$ to $\en$. Having obtained the data, we compute the posterior probability of the true model using the set of priors (\ref{4}) for several mixtures on $g$ as indicated above. There are two issues to be mentioned here. Firstly, for calculation of the marginal, one needs to calculate the integral in (\ref{8}), which is not of closed form for all the mixtures. We use numerical integration (available in R software) to calculate this integral for all the mixtures.
Secondly, since $p$ is large, calculation of the posterior probability for the candidate models $\alpha \in \scrptA$ becomes quite infeasible. This is because calculation of posterior probability of any model requires the marginal densities $\left(m_{\alpha}(\yn)\right)$ for all the $2^p$ candidate models in $\scrptA$. Therefore, we use Markov Chain Monte Carlo simulation techniques to approximate the posterior probabilities, where computation of marginal densities can be restricted only to the models visited by the chain. We have used the Gibbs sampling algorithm, to simulate from the relevant Markov chain. The sampling scheme and the method of computation of posterior probabilities are completely described in \citet[Section 3.5]{CGM_2001}. We have generated a Markov chain of length 10000 of which the first 5000 have been used as burn-in. Similar simulation is also used in Scheme 2 of \cite{MSC_2013}.

For each combination of $(n, p)$ and each of the mixtures of $g$, we repeat the above for 100 times fixing the chosen values of $\xi_i$'s, $p(\alpha_c)$, $\beta_j$'s and $\sigma$. The mean and mean square error 
of the posterior probabilities of the true model are presented in Table \ref{table:1}. In the table, the mean posterior probabilities have been shown, keeping the mean squared errors (m.s.e.) in brackets.

From Table \ref{table:1}, it is evident that performance of the proposed mixture is better than the other mixtures. The differences in performances of the proposed priors from others increase with $n$. Among the two choices of $\nu$, the choice $\nu=p$ performs better than $\nu=1$. Among the other priors, robust prior has the best performance.  
For example, when $t_{(3)}$ is considered with $n=150$ and $p=50$, for the proposed prior I (proposed prior II), the Markov chain visits the true model in around $37\%$ ($46\%$) cases, whereas for the robust prior it visits the true model only in $6\%$ cases. A similar phenomenon is observed in other cases as well.

\begin{figure}
  \caption{\footnotesize Mean and M.S.E. of Posterior Probabilities of the True Model in Scheme 1}
 \includegraphics[height=3.2 in,width=5.3 in]{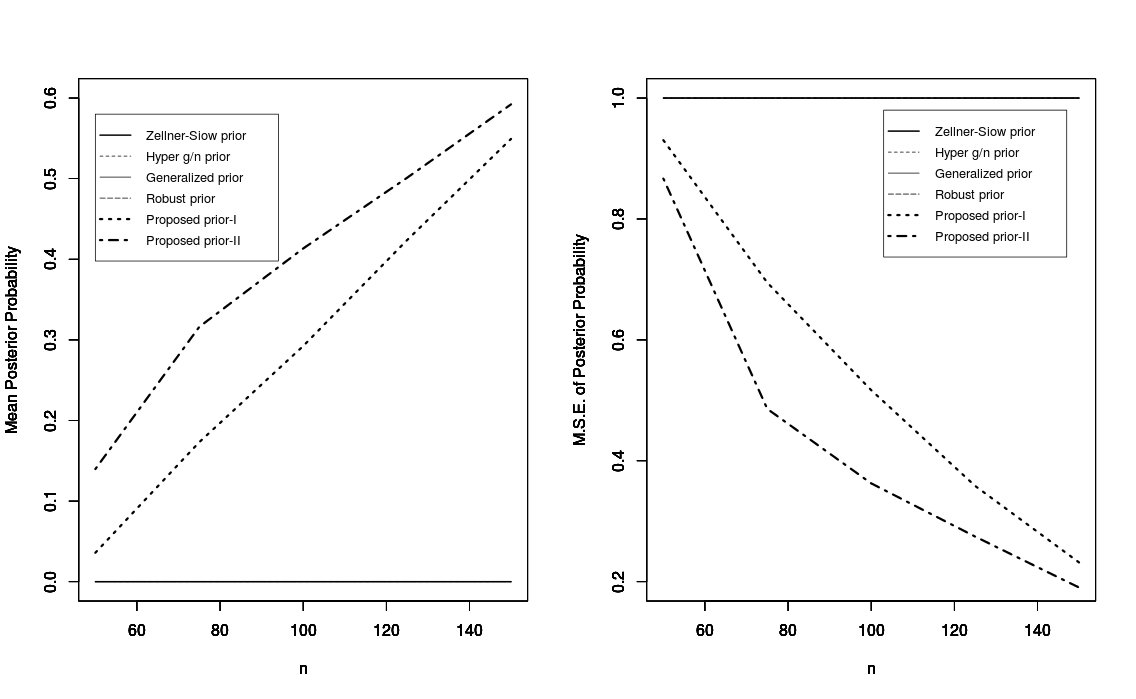}
 \label{fig3a}
\end{figure}

Next, we consider two sparse situations. We first consider the case where the null model is true ({\it Scheme 1}). We take $p=30$ and $\beta_0=5$  and assume $\en$ to follow a normal distribution. The rest of the simulation scheme is as described above. The mean and m.s.e. of 100 replicates of posterior probabilities of different mixtures on $g$ are shown in Figure \ref{fig3a}.

\begin{figure}
\centering
\caption{\footnotesize Mean and M.S.E. of Posterior Probabilities of the Sparse Model in Scheme 2}
\includegraphics[height=3.2 in,width=5.3 in]{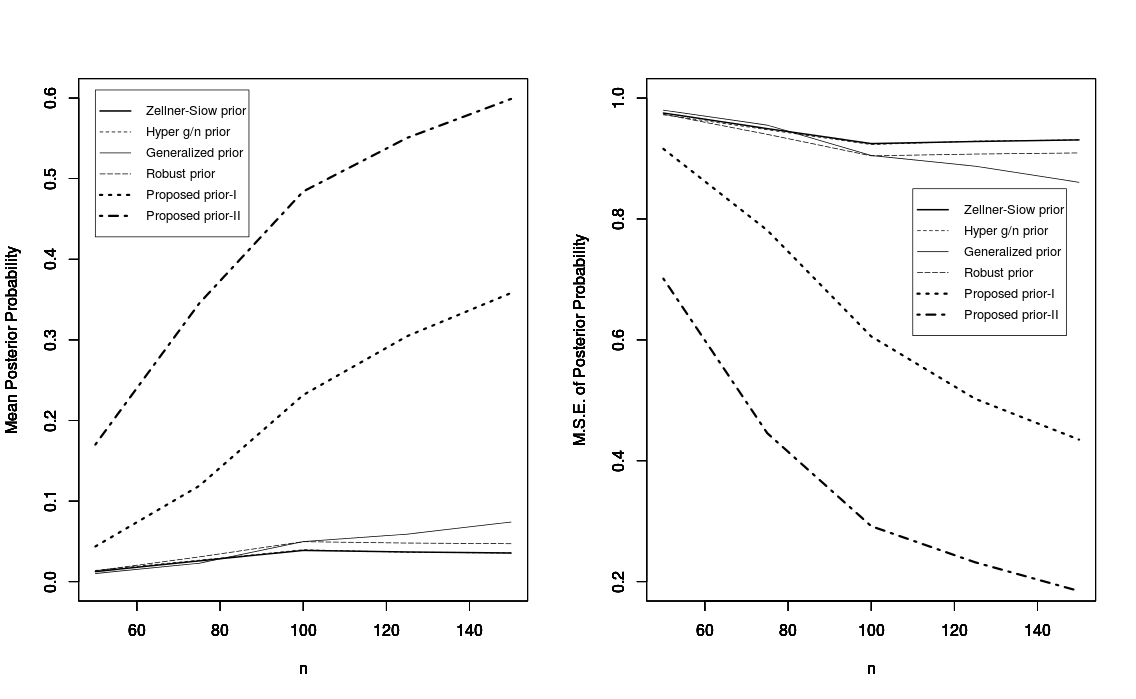}
 \label{fig3b}
\end{figure}

Lastly, we consider the interesting situation where sparsity is present in the simulation scheme ({\it Scheme 2}) in the sense that a set of regression coefficients in the true model is negligible, even if not exactly zero. Here, it is desirable to select the parsimonious model (say $\alpha_s$) that includes all regressors with significant regression parameters, rather than the true model ($\alpha_c$). We take $p=30$, $\pnac=15$, $\sigma^2=1$. The error distribution considered is normal. The regression parameters of the true model is as follows: $\beta_i=i+1$ for $i=0,1,\ldots , 4$, $0<|\beta_i|<0.008$ for $i=5,6,\ldots 15$ and $\beta_i=0$ for $i>15$. For different $n$, we plot the mean and m.s.e. of posterior probabilities of the sparse model, $\alpha_s=\{1,2,3,4\}$ in Figure \ref{fig3b}.     

 From both the figures, it can be seen that the performances of the proposed priors are distinctly better than that of the other priors. When the null model is true (see Figure \ref{fig3a}), the posterior probabilities of all other priors are less than $0.00002$ (therefore, the corresponding lines are not visually distinguishable in the figure), whereas the proposed prior for $\nu=p$ achieves an average posterior probability of $0.59$, when $n=150$. In {\it Scheme 2} (see Figure \ref{fig3b}), the other priors perform relatively better than {\it Scheme 1}. The highest average posterior probability among all the other priors is achieved by the generalized $g$-prior for $n=150$ and is less than $0.075$, whereas the average posterior probabilities achieved by the proposed priors for $n=150$ are $0.358$ for $\nu=1$, and $0.599$ for $\nu=p$.
   \section{Concluding Remarks}
In this paper, we propose a mixture of $g$-priors suitable for the case when $p$ grows with $n$. The resulting marginal has an approximation with a closed form expression which makes implementation simple. We investigate the performance of the proposed prior by deriving consistency properties under different setups. We also compare its performance with that of several other mixtures using simulation results under different simulation schemes which demonstrates its nobility.
Theoretically as well as in simulations, superiority of the performance of the proposed prior has also been shown under sparse situations.

The prior for $\Ba$ arising from this mixture has a very thick tail which is recommended by Jeffreys (\cite{Jeffreys_1961}). Further, this structure of priors (\ref{4}) has the properties like {\it predictive matching} and {\it group invariance} as described in \cite{Bayari_2012} (see Results 2-4 of \cite{Bayari_2012} in this context). The authors have explicitly justified the adoption of the form (\ref{4}) in a broader context.

Finally, it may be mentioned that we have studied the performance of the proposed mixture for $\nu=1$ and $\nu=p$.
The performance of the mixture with $\nu=p$ is better than the other in the light of all the properties considered in this paper except {\it information consistency}. The prior with $\nu=p$ fails to be information consistent when $n\leq 2p$. In practice, when $n>2p$ one can conveniently use the prior with $\nu=p$.  
 
\appendix
\section{}
In this section, we present the proofs of most of the main results stated in this paper. Many of the statements in the following proofs hold with probability tending to 1 as $n\rightarrow\infty$, although this will not be always mentioned. Throughout this section we will assume that $Var(\en)=\sigma^2 I$, where $\sigma^2>0$ is unknown.
\subsection{Auxiliary Results}
 We first state three lemmas which will help in proving our main results. The proofs of these lemmas are given in the supplementary file.

\begin{landscape}
\begin{small}
 \begin{table}[htbp] 
 \caption{Average and mean squared error of posterior probability of the {\it true} model.}
 \label{table:1}       
\begin{center}
\renewcommand{\arraystretch}{1.2}
\begin{tabular}{p{1 cm}llp{2.4 cm}p{2.4 cm}p{2.4 cm}p{2.4 cm}p{2.4 cm}p{2.4 cm}}
 \hline\noalign{\smallskip}
error  & \multicolumn{2}{c}{Priors} & Zeller-Siow & Hyper-$g/n$ & Generalized $g$ & Robust & Proposed I & Proposed II  \\ 
density & $p+1$ & $n$  & mean (m.s.e.)&mean (m.s.e.) &mean (m.s.e.) &mean (m.s.e.) &mean (m.s.e.) &mean (m.s.e.) \\
 \noalign{\smallskip}\hline\noalign{\smallskip}
  &$30$ & 50  & 0.0948 (0.8294) & 0.0949 (0.8293) & 0.0956 (0.8282) & 0.0981 (0.8242) & 0.1251 (0.7836) & 0.1598 (0.7356) \\
      & & 100 & 0.2916 (0.5095) & 0.2919 (0.5091) & 0.2926 (0.5081) & 0.3056 (0.4903) & 0.4598 (0.3131) & 0.6163 (0.1741) \\
Normal& & 150 & 0.4932 (0.2954) & 0.4962 (0.2922) & 0.3758 (0.4192) & 0.5358 (0.2515) & 0.6834 (0.1393) & 0.9055 (0.0174) \\
   &$50$ & 50 & 0.0000 (0.9999) & 0.0000 (0.9999) & 0.0000 (0.9999) & 0.0000 (0.9999) & 0.0006 (0.9999) & 0.0008 (0.9984) \\
      & & 100 & 0.0018 (0.9964) & 0.0017 (0.9965) & 0.0018 (0.9965) & 0.0019 (0.9963) & 0.0255 (0.9525) & 0.0414 (0.9265) \\
      & & 150 & 0.0474 (0.9086) & 0.0475 (0.9083) & 0.0477 (0.2749) & 0.0525 (0.8991) & 0.2581 (0.4382) & 0.4820 (0.3072) \\
  &$30$ & 50     & 0.1611 (0.7156) & 0.1603 (0.7167) & 0.1614 (0.7151) & 0.1669 (0.7067) & 0.2155 (0.6387) &  0.2767 (0.5547) \\
      & & 100    & 0.1937 (0.6620) & 0.1933 (0.6626) & 0.1944 (0.6608) & 0.2028 (0.6484) & 0.3374 (0.4764) &  0.4369 (0.3740) \\
Laplace& & 150   & 0.2573 (0.5640) & 0.2575 (0.5636) & 0.2578 (0.5632) & 0.2734 (0.5419) & 0.4822 (0.3073) &  0.6144 (0.2095) \\
  &$50$ & 50  & 0.0011 (0.9978) & 0.0011 (0.9978) & 0.0011 (0.9978) & 0.0011 (0.9978) & 0.0013 (0.9974) &  0.0068 (0.9866) \\
      & & 100 & 0.0052 (0.9896) & 0.0053 (0.9895) & 0.0053 (0.9895) & 0.0055 (0.9890) & 0.0585 (0.8941) &  0.0816 (0.8584) \\
      & & 150 & 0.0069 (0.9863) & 0.0070 (0.9861) & 0.0070 (0.9862) & 0.0075 (0.9851) & 0.1166 (0.8024) &  0.1563 (0.7495) \\
&$30$ & 50   & 0.1059 (0.8103) & 0.1060 (0.8103) & 0.1065 (0.8094) & 0.1091 (0.8053) & 0.1673 (0.7252) &  0.1921 (0.6918) \\
      & & 100& 0.3191 (0.4700) & 0.3197 (0.4693) & 0.3204 (0.4683) & 0.3342 (0.4497) & 0.5003 (0.2674) &  0.6593 (0.1332) \\
$t_{(3)}$&&150 & 0.5419 (0.2103) & 0.5460 (0.2488) & 0.4292 (0.3629) & 0.5834 (0.2141) & 0.7100 (0.1190) &  0.8965 (0.0274) \\
  &$50$ & 50  & 0.0004 (0.9993) & 0.0004 (0.9993) & 0.0004 (0.9992) & 0.0004 (0.9992) & 0.0307 (0.9444) & 0.0254 (0.9522) \\
      & & 100 & 0.0080 (0.9840) & 0.0081 (0.9840) & 0.0080 (0.9840) & 0.0084 (0.9834) & 0.0713 (0.8765) & 0.0911 (0.8468) \\
      & & 150 & 0.0562 (0.8924) & 0.0562 (0.8924) & 0.0557 (0.8931) & 0.0615 (0.8826) & 0.3744 (0.4279) & 0.4653 (0.3316) \\
 \noalign{\smallskip}\hline
 \end{tabular}
 \end{center}
 \end{table}
 \end{small}
 \end{landscape}

\begin{lemma}\label{lm_2}
If $\yn=\mun+\en$ with $\mun$ satisfying assumption (A.1) and $\en$ satisfying assumption (A.2), then the following results hold as $n\rightarrow \infty$:

$(i)$ $ \bar{e}=\sum_{i=1}^n e_i /n \xrightarrow{p} 0$,

$(ii)$ $\sum_{i=1}^{n} \mu_i e_i /n \xrightarrow{p} 0$,

$(iii)$ $\sesq \xrightarrow{p}\sigma^2$ where $n\sesq=\sum _{i=1}^{n}\left(e_i - \bar{e}\right)^2$, 

$(iv)$ $\max_{\alpha\in\scrptA} \en^{\prime}\hatmat\en/n = O_p \left(p/n\right)$,

$(v)$ $\max_{\alpha\in\scrptA_2}\left|\mun^\prime(I-\hatmat)\en\right|/n=O_p\left(\sqrt{p/n}\right)$ and 

$(vi)$ $\en^\prime(I-\hatmat)\en/n\xrightarrow{p}\sigma^2$ uniformly in $\alpha\in\scrptA$.  
\end{lemma}

\begin{lemma}\label{lm_1} Let $\rasq$ is as in (\ref{8}). Under assumptions (A.1) and (A.2), $(1-\rasq)>\sigma^2/\left( 2M+ 4 \sigma^2\right)$ with probability tending to 1 uniformly in $\alpha\in\scrptA$, where $M$ is as in assumption (A.1). 
\end{lemma}
\begin{lemma}
 \label{lm_3}
Under the setup of Theorem \ref{thm_3}, for any fixed $R>0$, with probability tending to one $$\max_{\alpha\in\scrptA_1^{*}}\frac{\en^{\prime}(\hatmat-\hatmatAc)\en}{\sigma^2(\pna-\pnac)} \leq R \log p.$$
\end{lemma}

\subsection{Proof of Result \ref{res_1}}
\noindent Using (\ref{5}) and (\ref{8}), we write,
\begin{eqnarray}
 m_{\alpha}(\yn)=\mathcal{C}_{1,y,n} \mathcal{I}  \label{a.1}
\end{eqnarray}
where,
\begin{eqnarray}
 \mathcal{C}_{1,y,n}=\frac{\Gamma{(n-1)/2}}{\Gamma{(\nu/2)}{\pi}^{(n-1)/2}\sqrt{n}}\left( \frac{{\tau}^2\nu}{2}\right)^{\nu/2}\left(\sysq\right)^{-(n-1)/2}, \label{a.2}
\end{eqnarray}
and
\begin{equation}
 \mathcal{I}=\int_{0}^{\infty} e^{-{\tau^2\nu/(2g)}} g^{-(1+\nu/2)} (1+g)^{(n-\pna-1)/2}\left\{ 1+(1-\rasq)g\right\}^{-(n-1)/2} dg. \label{a.3}
\end{equation}
We first evaluate $\mathcal{I}$. After making a transformation $w=\tau^2\nu/(2g)$, we observe that,
\begin{equation}
 \mathcal{I}=\mathcal{C}_{2,y,n}\int_{0}^{\infty} e^{-w} w^{\nu/2-1} \left(1+\frac{\tau^2\nu}{2w}\right)^{(n-\pna-1)/2} \left\{1+(1-\rasq)\frac{\tau^2\nu}{2w} \right\}^{-(n-1)/2} dw, \label{a.8}
\end{equation}
where $\mathcal{C}_{2,y,n} =\left(\tau^2\nu/2\right)^{-\nu/2} $. Next we use the fact that, for any $w>0$,
$$ \left\{1+(1-\rasq)\tau^2\nu/(2w) \right\}^{-(n-1)/2} < \left\{(1-\rasq)\tau^2\nu/(2w) \right\}^{-(n-1)/2}.$$
Use of this inequality along with multiplication and division by $\left(\tau^2\nu/(2w)\right)^{(n-\pna-1)/2}$  in R.H.S. of (\ref{a.8}) gives,
\begin{eqnarray}
 \mathcal{I}\leq \mathcal{C}_{3,y,n} \int_{0}^{\infty} e^{-w}  w^{(\pna+\nu)/2-1} \left(1+\frac{2w}{\tau^2\nu}\right)^{(n-\pna-1)/2} dw, \label{a.4}
\end{eqnarray}
where $\mathcal{C}_{3,y,n} =\mathcal{C}_{2,y,n} (\tau^2\nu/2)^{-\pna/2}(1-\rasq)^{-(n-1)/2}$. The quantity $(n-\pna-1)/2$ can either be an integer or a mixed fraction. We first deal with the case when $(n-\pna-1)/2$ is an integer. We expand last term in (\ref{a.4}) in binomial expansion as follows, 
\begin{align}
\left(1+\frac{2w}{\tau^2\nu}\right)^{(n-\pna-1)/2} = \left( 1+ \frac{(n-\pna-1) w}{\tau^2\nu}+
\ldots + \left(\frac{2w}{\tau^2\nu}\right)^{(n-\pna-1)/2} \right). \label{a.5} 
\end{align}
 For a fixed $n$ the sum in (\ref{a.5}) is finite and hence the integration in (\ref{a.4}) and the summation can be interchanged. It can be easily seen that after interchange has taken place, each integration under the summation forms a gamma integral of appropriate order. It then follows that,
 \begin{eqnarray}
  \mathcal{I}&\leq&  \mathcal{C}_{3,y,n}  \left\{\Gamma{\left(\frac{\pna+\nu}{2}\right)}+ \frac{(n-\pna-1)}{\tau^2\nu}\Gamma{\left(\frac{\pna+\nu}{2}+1\right)} + \right.\nonumber \\
&&~\left. \hspace{1.6 in}
\ldots + \left( \frac{2}{\tau^2\nu}\right)^{(n-\pna-1)/2} \Gamma{\left(\frac{n+\nu-1}{2}\right)} \right\}, \nonumber\\
&=& \mathcal{C}_{3,y,n} \Gamma{\left(\frac{\pna+\nu}{2}\right)} \left\{ 1+ \frac{(n-\pna-1) (\pna+\nu)}{2\tau^2\nu}+ \right.\nonumber\\
&&~\left. \hspace{.9 in}
\ldots +\frac{(\pna+\nu)(\pna+\nu+2)+\ldots +(n+\nu-3)}{\left(\tau^2\nu\right)^{(n-\pna-1)/2}} \right\}. \nonumber
 \end{eqnarray}
For $\tau=n$, the above expression is,
\begin{eqnarray}
 &<& \mathcal{C}_{3,y,n} \Gamma{\left(\frac{\pna+\nu}{2}\right)} \left\{ 1+ \frac{n (\pna+\nu)}{2n^2\nu}+\frac{n^2(\pna+\nu)(\pna+\nu+2)}{2! \left(2 n^2\nu\right)^2}+\right.\nonumber\\
&&~\left. \hspace{1.8 in} \ldots +\mbox{~~~~~upto $(n-\pna-1)/2^{th}$ term } \right\} \nonumber \\
 &<& \mathcal{C}_{3,y,n} \Gamma{\left(\frac{\pna+\nu}{2}\right)} \left\{ 1+ \frac{\pna+\nu}{2n\nu}+\left(\frac{\pna+\nu}{2n\nu}\right)^2+ \right.\nonumber\\
&&~\left. \hspace{2.4 in}\ldots +\left(\frac{\pna+\nu}{2n\nu}\right)^{(n-\pna-1)/2} \right\} \nonumber
\end{eqnarray}
using the fact that if  $a>b>0$ then $a/b>(a+1)/(b+1)$.
The bracketed portion of the last expression is a G.P. series with positive terms. We add the terms with higher power and make it an infinite G.P. series and also replace $\pna$ by $p$ ($\pna\leq p$) to make the series free of $\alpha$. The resultant term is as follows, 
\begin{eqnarray}
 &<& \mathcal{C}_{3,y,n} \Gamma{\left(\frac{\pna+\nu}{2}\right)} \left\{ 1 + \frac{p+\nu}{2n\nu}+\left(\frac{p+\nu}{2n\nu} \right)^{2}+\ldots\right\} \label{a.6}
\end{eqnarray}
From (\ref{a.6}) it is clear that 
\begin{align}
 \mathcal{I}\leq\left\{ \begin{array}{ll} \mathcal{C}_{3,y,n} \Gamma{\left(\frac{\pna+\nu}{2}\right)} \left(1+O\left(\frac{1}{n^{1-b}}\right) \right) \quad \quad \mbox{if $\nu=1$ and $p=n^b$,}\\
 \mathcal{C}_{3,y,n} \Gamma{\left(\frac{\pna+\nu}{2}\right)} \left(1+O\left(\frac{1}{n}\right) \right)\quad \quad \quad\mbox{if $\nu=p$ and $p=n^b$.}\\
 \end{array} \right. \label{a.7}
\end{align}
We now consider the case when $(n-\pna-1)/2$ is not an integer. Certainly then $(n-\pna)/2$ is an integer. For any $w>0$, we write the last term in (\ref{a.4}) as
\begin{eqnarray*}
\left(1+(2w)/(\tau^2\nu)\right)^{(n-\pna -1)/2}
\leq \left(1+(2w)/(\tau^2\nu)\right)^{(n-\pna)/2}.
\end{eqnarray*}
Using the above inequality, we proceed as before in (\ref{a.5}), (\ref{a.6}) and get the same bound as in (\ref{a.7}).

Next, we assign a bound on $\mathcal{I}$ from other direction and show that the difference between the two bound is small. For this, we move back to (\ref{a.8}) and use the inequality $\left(1+(\tau^2\nu)/(2w)\right) > (\tau^2\nu)/(2w)$ along with a multiplication and division by the factor $\left((1-\rasq)(\tau^2\nu)/(2w)\right)^{(n-1)/2}$ in the integrand of (\ref{a.8}). The resultant integral is as follows,
\begin{eqnarray*}
 \mathcal{I} &\geq& \mathcal{C}_{3,y,n} \int_{0}^{\infty} w^{(\nu+\pna)/2-1} e^{-w} \left(1-\frac{1}{1+(1-\rasq)\tau^2\nu/(2w)} \right)^{(n-1)/2} dw, \\
  &\geq& \mathcal{C}_{3,y,n} \int_{0}^{\infty} w^{(\nu+\pna)/2-1} e^{-w} \left(1-\frac{2w}{(1-\rasq)\tau^2\nu} \right)^{(n-1)/2} dw,
\end{eqnarray*}
where $\mathcal{C}_{3,y,n}$ is the same as in (\ref{a.7}). As before, here also we deal separately two cases when $(n-1)/2$ is an integer and when it is not. First consider the case when it is an integer. Expanding as before we get,
\begin{eqnarray*}
 \left(1-\frac{2w}{(1-\rasq)\tau^2\nu} \right)^{(n-1)/2} \hspace{4.1 in}\\
 = 1-\frac{(n-1)w}{\tau^2\nu(1-\rasq)} +\frac{(n-1)(n-3)w^2}{2!\tau^4\nu^2(1-\rasq)^2}- \ldots +\left(\frac{-2w}{(1-\rasq)\tau^2\nu}\right)^{(n-1)/2}  \hspace{1.2 in}\\
   \geq 1-\frac{w}{n\nu(1-\rasq)} -\frac{1}{2!}\left(\frac{w}{(1-\rasq)n\nu}\right)^{2}- \ldots -\frac{1}{((n-1)/2)!}\left(\frac{w}{(1-\rasq)n\nu}\right)^{(n-1)/2}, \quad \quad
\end{eqnarray*}
putting $\tau=n$ and replacing $(n-1), ~(n-3),$ etc. by $n$.
As before, we interchange the summation and integration and resultant term is as follows
\begin{eqnarray}
 \mathcal{I} \geq \mathcal{C}_{3,y,n}\Gamma{\left(\frac{\pna+\nu}{2}\right)}\left\{1-\frac{(\pna+\nu)}{2n\nu(1-\rasq)} -\frac{(\pna+\nu)(\pna+\nu+2)}{2!\left(2n\nu(1-\rasq)\right)^2}- \ldots \right\}\hspace{.2 in} \nonumber \\
 \hspace{.2 in}\geq \mathcal{C}_{3,y,n}\Gamma{\left(\frac{\pna+\nu}{2}\right)}\left\{1-\frac{(\pna+\nu)}{2n\nu(1-\rasq)} -\left(\frac{(\pna+\nu)}{2n\nu(1-\rasq)}\right)^2- \right.\hspace{1 in} \nonumber\\
 \left.  \ldots - \left(\frac{(\pna+\nu)}{2n\nu(1-\rasq)}\right)^{(n-1)/2}  -\ldots \right\}\hspace{.4 in} \nonumber\\
 \hspace{.1 in} \geq \mathcal{C}_{3,y,n}\Gamma{\left(\frac{\pna+\nu}{2}\right)}\left\{1-\frac{(p+\nu)}{2n\nu(1-\rasq)} -\left(\frac{(p+\nu)}{2n\nu(1-\rasq)}\right)^2 -\ldots~\right\}.\hspace{.6 in} \nonumber
\end{eqnarray}
From Lemma \ref{lm_1}, we know that for all $\alpha\in\scrptA$, $\left(1-\rasq\right)$ has a fixed positive lower bound with probability tending to 1. Using the lemma we get the following result
\begin{align}
 \mathcal{I}\geq\left\{ \begin{array}{ll} \mathcal{C}_{3,y,n} \displaystyle\Gamma{\left(\frac{p+\nu}{2}\right)} \left(1+O_{p}\left(\frac{1}{n^{1-b}}\right) \right) \quad  \mbox{if $\nu=1$ and $p=n^b$,}\\
 \mathcal{C}_{3,y,n} \displaystyle\Gamma{\left(\frac{p+\nu}{2}\right)} \left(1+O_{p}\left(\frac{1}{n}\right) \right) \quad \quad\mbox{if $\nu=p$ and $p=n^b$.}\\
 \end{array} \right. \label{a.10}
\end{align}
Hence, from (\ref{a.1}), (\ref{a.7}) and (\ref{a.10}), the result follows. \qed

\subsection{Proof of Theorem \ref{thm_1}}\label{Sect1}
We need to show
\begin{equation}
 \sum_{\alpha\in \scrptA_i}\frac{p(\Ma)}{p(\Mac)} \frac{m_{\alpha}(\yn)}{m_{\alpha_c}(\yn)} \xrightarrow{p} 0, ~~~~\mbox{for } i=1,2.\label{a.12}
\end{equation}
We proof separately for the cases $M_{\alpha_c}\neq M_{N}$ and $M_{\alpha_c}=M_{N}$.

\noindent {\it Case I. $M_{\alpha_c}\neq M_{N}$.}
We first consider (\ref{a.12}) for $i=2$.

From Result \ref{res_1} we have 
\begin{eqnarray}
 \frac{m_{\alpha}(\yn)}{m_{\alpha_c}(\yn)} &\leq& \left(\frac{n^2\nu}{2}\right)^{-(\pna-\pnac)/2} \left(\frac{1-R_{\alpha_c}^2}{1-\rasq} \right)^{(n-1)/2} \nonumber \\
 && \hspace{1 in} \frac{\Gamma{\left\{(\nu+\pna)/2\right\}}}{\Gamma{\left\{(\nu+\pnac)/2\right\}}} 
\frac{\left\{1+p O(1)/(\nu n)\right\}}{\left\{1+p O_p(1)/(\nu n) \right\}}, \label{a.13}
\end{eqnarray}
where the terms $O(1)$ and $O_p(1)$ are free of $\alpha$.
We consider the terms of R.H.S. of (\ref{a.13}) one by one. First we consider the second term. We have for $\alpha\in\scrptA_2$,
\begin{eqnarray}
 &&\left(\frac{1-\rasq}{1-R_{\alpha_c}^2} \right)\nonumber \\
&& =\frac{\en^{\prime}\en +\mun^{\prime}\left(I-\hatmat \right)\mun+2\mun^{\prime}\left(I-\hatmat \right)\en-\en^{\prime}\hatmat\en}{\en^{\prime}\left(I-\hatmatAc \right)\en},\quad \quad \quad\quad \quad \quad \nonumber\\
\label{a.15} &&\geq  1+\frac{1}{\en^{\prime}\en/n}\left(\displaystyle\min_{\alpha\in\scrptA_2}\frac{\mun^{\prime}\left(I-\hatmat \right)\mun}{n}-2\displaystyle\max_{\alpha\in\scrptA_2} \frac{\left|\mun^{\prime}\left(I-\hatmat \right)\en\right|}{n} \right. \nonumber\\ 
&& \left. \hspace{3.2 in} -\displaystyle\max_{\alpha\in\scrptA_2}\frac{\en^{\prime} \hatmat\en}{n}  \right).\nonumber
\end{eqnarray}
By assumption (A.3) and from parts $(iv)$ and $(v)$ of Lemma \ref{lm_2}, we have 
$$\frac{1}{n}\left\{ \min_{\alpha\in\scrptA_2}\mun^{\prime}\left(I-\hatmat \right)\mun- \max_{\alpha\in\scrptA_2} \left( 2\left|\mun^{\prime} \left(I-\hatmat \right)\en\right| + \en^{\prime} \hatmat\en \right) \right\} >\frac{\delta_0}{n^s}$$
for some $\delta_0>0$, with probability tending to 1 provided $s<(1-b)/2$.
Since $\en^{\prime}\en/n\xrightarrow{p}\sigma^2$,
\begin{equation}
 \displaystyle\max_{\alpha\in\scrptA_2}\left(\frac{1-R_{\alpha_c}^2}{1-\rasq} \right)^{(n-1)/2} \leq \left( 1+\frac{\delta_1}{n^s} \right)^{-(n-1)/2} \hspace{.2 in}\label{a.34}
\end{equation}
for some $\delta_1>0$.
 
 To evaluate third term of R.H.S. of (\ref{a.13}), we make use of the result $\{x/(x+s)\}^s\leq \Gamma{(x+s)}/\left(x^s \Gamma{x}\right) \leq 1$ for $0<s<1$ and $x>0$ from \cite{Gm_bound}. It can be shown that,
\begin{eqnarray}
 \frac{\Gamma{\left\{(\nu+\pna)/2\right\}}}{\Gamma{\left\{(\nu+\pnac)/2\right\}}}\leq\left(\frac{\nu+p}{2}\right)^{\left|\pna-\pnac\right|/2} . \label{a.16}
\end{eqnarray}
Hence from (\ref{a.13}), (\ref{a.34}) and (\ref{a.16}),
\begin{eqnarray*}
 &&\max_{\alpha\in\scrptA_2} \frac{m_{\alpha}(\yn)}{m_{\alpha_c}(\yn)} \\
 && \leq \left(\frac{n^2\nu}{2}\right)^{(p-\pnac)/2} \left(1+\frac{\delta_1}{n^s}\right)^{-(n-1)/2} \left(\frac{\nu+p}{2}\right)^{p/2} \frac{\left\{1+p ~O(1)/(\nu n)\right\}}{\left\{1+p ~O_p(1)/(\nu n) \right\}}. 
\end{eqnarray*}
 Using assumption (A.4) we also get an upper bound of the ratio of prior probabilities of the models. Therefore,
\begin{align}
 \label{a.22} \sum_{\alpha\in\scrptA_2} \frac{p(\Ma)}{p(\Mac)}\frac{m_{\alpha}(\yn)}{m_{\alpha_c}(\yn)} \leq C^{\prime} 2^{p} \left(\frac{n^2\nu}{2}\right)^{(p-\pnac)/2} \left(1+\frac{\delta_1}{n^s}\right)^{-(n-1)/2} \left(\frac{\nu+p}{2}\right)^{p/2}, 
\end{align}
for some constant $C^{\prime}$.
It is easy to check that the above quantity goes to $0$ as $n\rightarrow\infty$ when $s<(1-b)/2$.
\vskip5pt

Next we prove (\ref{a.12}) for $i=1$. We recall (\ref{a.13}) and consider each term of R.H.S.. We have, for any $\alpha\in\scrptA_1$,
\begin{eqnarray}
 \left(\frac{1-R_{\alpha_c}^2}{1-\rasq} \right)^{(n-1)/2} &=&\left\{ \frac{\en^{\prime}\left(I-\hatmatAc \right)\en}{\en^{\prime}\left(I-\hatmat \right)\en}\right\}^{(n-1)/2} \nonumber\\
&=&  \left\{1-\frac{\en^{\prime}\left(\hatmat-\hatmatAc \right)\en}{\en^{\prime}\left(I-\hatmatAc \right)\en} \right\}^{-(n-1)/2}  \nonumber
\end{eqnarray}
We now use Lemma 2 of \cite{MSC_2013} which states that for any $R>2$, with probability tending to 1,
\begin{equation}
  \max_{\alpha\in\scrptA_1} \frac{\en^{\prime}(\hatmat-\hatmatAc)\en }{ \sigma^2 (\pna-\pnac)} \leq R \log p. \label{a.23} 
\end{equation}
Again from part $(vi)$ of Lemma \ref{lm_2}, $\en^{\prime}\left(I-\hatmatAc \right)\en>n\sigma^2 (1-\epsilon)$ with probability tending to 1, for any $\epsilon>0$. These facts will imply that for any $R>2$, with probability tending to 1 uniformly in $\alpha\in\scrptA_1$,
\begin{eqnarray}
\left(\frac{1-R_{\alpha_c}^2}{1-\rasq} \right)^{(n-1)/2}&\leq&  
\left\{1-\frac{\sigma^2(\pna-\pnac)R\log p}{n\sigma^2(1-\epsilon)} \right\}^{-(n-1)/2} \nonumber\\
&\leq& \exp \left\{\frac{(\pna-\pnac)R\log p}{2(1-\epsilon)^2} \right\}  \nonumber\\
&&(\mbox{for any $0<z<\epsilon<1$,  $(1-z)>\exp\left\{-z/(1-\epsilon)\right\}$})\nonumber\\
&=& {\large (p)^{R(\pna-\pnac)/2(1-\epsilon)^2}}. \label{a.17}
\end{eqnarray}
Combining (\ref{a.13}), (\ref{a.16}), (\ref{a.17}) and using assumption (A.4) we have,
\begin{eqnarray}
 \sum_{\alpha\in\scrptA_1} \frac{p(\Ma)}{p(\Mac)} \frac{m_{\alpha}(\yn)}{m_{\alpha_c}(\yn)} &\leq& C \sum_{\alpha\in\scrptA_1}\left(\frac{(\nu+p){p}^{R/(1-\epsilon)^2}}{n^2\nu} \right)^{(\pna-\pnac)/2} \nonumber\\
  &\leq& C \sum_{q=1}^{p-\pnac}\binom{p-\pnac}{q}\left(\frac{\sqrt{\nu+p}~~{p}^{R/\{2(1-\epsilon)^2\}}}{n\sqrt{\nu}} \right)^{q} \nonumber\\
 &\leq& C \left\{ \left(1+\frac{\sqrt{\nu+p}~~{p}^{R/\{2(1-\epsilon)^2\}}}{n\sqrt{\nu}} \right)^{(p-\pnac)}-1 \right\} \label{a.18}
\end{eqnarray}
The above expression converges to $0$ as $n\rightarrow\infty$, if the first term in the curly braces converges to $1$. If $\nu=1$ and $p=O(n^b)$ then the first term is less than $\left( 1+C^{\prime}~n^{-[1-(R+1)b/\{2(1-\epsilon)^2\}]}\right)^{kn^b}$ for some positive constants $C^{\prime}$ and $k$, any $R>2$ and any $\epsilon>0$. Also if $\nu=p$ then this term is less than $\left( 1+C^{\prime\prime}~n^{-[1-Rb/\{2(1-\epsilon)^2\}]}\right)^{kn^b}$ for some positive constants $C^{\prime\prime}$ and $k$, any $R>2$ and any $\epsilon>0$. Letting $R\downarrow 2$ and $\epsilon\downarrow0$, we conclude that the last expression in (\ref{a.18}) converges to $0$ if $b<2/5$ when $\nu=1$ and if $b<1/2$ when $\nu=p$. 
\vskip5pt

{\it Case II. $M_{\alpha_c}=M_{N}$.}
When the null model is true, the Bayes factor of any model with respect to the null model is given by,
\begin{align}
  \label{a.20} \frac{m_{\alpha}(\yn)}{m_{N}(\yn)}=\int_{0}^{\infty} (1+g)^{(n-\pna-1)/2} \left\{ 1+g(1-\rasq)\right\}^{-(n-1)/2} \pi(g) dg,
\end{align}
where $\pi(g)$ is as in (\ref{5}). 
Now
\begin{eqnarray*}
 \left[\frac{1+g}{1+(1-\rasq)g}\right]^{(n-1)/2}
 &=&\exp\left[\left(\frac{n-1}{2}\right)\left\{\ln(1+g)-\ln(1+(1-\rasq)g) \right\} \right]  \\
 &=&\exp\left[\left(\frac{n-1}{2}\right)\left\{\ln(1+g)-\ln(1+g)+\frac{\rasq g}{1+g^*} \right\} \right] \quad\\
 &&\quad \quad \quad \quad\mbox{ where $g^*\in\left[(1-\rasq)g,g\right]$}\\
&\leq&\exp\left[\left(\frac{n-1}{2}\right)\left\{\frac{\rasq}{(1-\rasq)+1/g}\right\} \right] \\
&\leq &\exp\left[\left(\frac{n-1}{2}\right)\left(\frac{\rasq}{1-\rasq}\right) \right] \\
&=&\exp\left[\left(\frac{n-1}{2}\right)\left(\frac{\en^{\prime}(\hatmat-\hatmatAc)\en}{\en^{\prime}(I-\hatmat)\en}\right) \right]
\end{eqnarray*}
since $n\sysq=\en^{\prime}\en - n\bar{e}^2$ and $n\bar{e}^2=\en^{\prime}\hatmatAc\en$ when the null model is true.

Next we use the facts that with probability tending to 1, $\en^{\prime}(I-\hatmat)\en >(n-1)\sigma^2 (1-\epsilon)$, for any $\epsilon >0$ and for any $R>2$, $\max_{\alpha\in\scrptA}\en^{\prime}\hatmat\en)/\sigma^2 <R\pna\ln(p)$. Combining these, we have for any $\epsilon>0$,
\begin{eqnarray}
\left[\frac{1+g}{1+(1-\rasq)g}\right]^{(n-1)/2} \leq p^{R\pna/{2(1-\epsilon)}} \nonumber
\end{eqnarray}
with probability tending to 1.

Hence from (\ref{a.20})  we have,
\begin{eqnarray}
 \frac{m_{\alpha}(\yn)}{m_{N}(\yn)} \leq p^{R\pna/{2(1-\epsilon)}}  \int_{0}^{\infty} (1+g)^{-\pna/2} \pi(g) dg. \label{a.21}
\end{eqnarray}
Now with $\pi(g)$ as given in (\ref{5})
\begin{eqnarray}
 I= \int_{0}^{\infty} (1+g)^{-\pna/2} \pi(g) dg &\leq& \frac{\left( \tau^2 \nu/2\right)^{\nu/2}}{\Gamma{(\nu/2)}} \int_0^{\infty} e^{-\tau^2 \nu/2g} g^{-(\pna+\nu)/2-1} dg \nonumber
\end{eqnarray}
by the fact that $(1+g)^{-1}<g^{-1}$. We then have
\begin{eqnarray*}
 I &=& \left\{\begin{array}{ll} \left(\tau^2/2\right)^{-\pna/2} \Gamma{\{(\pna+1)/2\}}/\Gamma{(1/2)} & \mbox{for } \nu =1,\\
\left(\tau^2p/2\right)^{-\pna/2} \Gamma{\{(\pna+p)/2\}}/\Gamma{(p/2)} & \mbox{for } \nu =p.\\
\end{array}\right.
 \end{eqnarray*}
 To evaluate these terms, we again use the results of \cite{Gm_bound} stated above. After little algebra one can show
\begin{eqnarray*}
I &<& \left\{\begin{array}{ll} \left(\tau^2/2\right)^{-\pna/2} \left(p/2\right)^{\pna/2} & \mbox{for } \nu =1,\\
\left(\tau^2p/2\right)^{-\pna/2} p^{\pna/2} & \mbox{for } \nu =p.\\
\end{array}\right.
 \end{eqnarray*}
 By assumption (A.4) and putting $\tau=n$, it follows from (\ref{a.21}) that
 \begin{eqnarray}
 &&\sum_{\alpha\in\left(\scrptA-\{\alpha_c\}\right)} \frac{p(\Ma)}{p(\Mac)} \frac{m_{\alpha}(\yn)}{m_{\alpha_c}(\yn)}\nonumber \\
 && ~~< \left\{\begin{array}{ll} C\displaystyle\sum_{\alpha\in\left(\scrptA-\{\alpha_c\}\right)} \left(p^{1+R/(1-\epsilon)}/n^2\right)^{\pna/2}/\sqrt{\pi} & \mbox{for } \nu =1,\nonumber\\
C\displaystyle\sum_{\alpha\in\left(\scrptA-\{\alpha_c\}\right)} \left(2p^{R/(1-\epsilon)}/n^2\right)^{\pna/2} & \mbox{for } \nu =p.\\
\end{array}\right.\\
\label{a.27} &&~~= \left\{\begin{array}{ll} C\left\{ \left(1+p^{\{1+R/(1-\epsilon)\}/2}/n\right)^{p}-1 \right\} & \mbox{for } \nu =1,\\
 C\left\{ \left(1+\sqrt{2}p^{R/\{2(1-\epsilon)\}}/n\right)^{p}-1 \right\}& \mbox{for } \nu =p,\\
\end{array}\right.
 \end{eqnarray}
 for any $R>2$ and any $\epsilon>0$.
 
As before, we let $R\downarrow 2$ and $\epsilon \downarrow 0$ and observe that the above quantity converges to $0$ when $p$ is of order $n^b$ if $b<2/5$ for $\nu=1$ and $b<1/2$ for $\nu=p$. \qed

\subsection{Proof of Theorem \ref{thm_3}} We proceed as in the proof of the Theorem \ref{thm_1}, and prove (\ref{a.12}) with $\scrptA_i$ replaced by its analog for nested models, $\scrptA_i^{*}$.
As before we consider separately the cases when the true model is null and when it is non-null.

{\it Case-I: $M_{\alpha_c} \neq M_N$.} First we consider $i=2$. Recall the proof of Theorem \ref{thm_1}. It can easily be seen that this part in Theorem \ref{thm_1} is proved for the model space $\scrptA_2$ and without using the assumption of normality. Since $\scrptA_2^{*}$ is a proper subset of $\scrptA_2$, here also the same proof holds.

Next consider (\ref{a.12}) with $\scrptA_i$ replaced by $\scrptA_i^{*}$ and $i=1$.
 Here we use Lemma \ref{lm_3} which is equivalent to (\ref{a.23}) when nested models are considered. This implies that (\ref{a.17}) also holds here. 
Combining (\ref{a.13}), (\ref{a.16}) and (\ref{a.17}) we have, for any $R>0$ and any $\epsilon>0$,

\begin{eqnarray}
 &&\sum_{\alpha\in\scrptA_1^{*}} \frac{p(\Ma)}{p(\Mac)}\frac{m_{\alpha}(\yn)}{m_{\alpha_c}(\yn)} \nonumber\\
 &&~\leq C\sum_{\alpha\in\scrptA_1}\left(\frac{(\nu+p){p}^{R/(1-\epsilon)^2}}{n^2\nu} \right)^{(\pna-\pnac)/2} \nonumber\\
 &&~\leq C \sum_{q=1}^{p-\pnac}\left(\frac{\sqrt{\nu+p}~~{p}^{R/2(1-\epsilon)^2}}{n\sqrt{\nu}} \right)^{q} \nonumber
 \end{eqnarray}
 \begin{eqnarray}
 &&~\leq C\left(\frac{\sqrt{\nu+p}~~{p}^{R/2(1-\epsilon)^2}}{n\sqrt{\nu}} \right) \frac{1}{\left(1-\sqrt{\nu+p}~~{p}^{R/2(1-\epsilon)^2}/n\sqrt{\nu} \right)}. \label{a.26}
\end{eqnarray}
For suitably chosen $R$ and $\epsilon$, it can be easily seen that (\ref{a.26}) converges to $0$ for any $0<b<1$.
\vskip5pt

{\it Case II. $M_{\alpha_c}=M_N$.} We proceed as in Theorem \ref{thm_1}. Observe that using Lemma \ref{lm_3} we obtain the inequality similar to (\ref{a.27}) with $\scrptA$ replaced by $\scrptA^{*}$. Thus we have
\begin{eqnarray*}
 &&\sum_{\alpha\in\left(\scrptA^{*}-\{\alpha_c\}\right)} \frac{p(\Ma)}{p(\Mac)} \frac{m_{\alpha}(\yn)}{m_{\alpha_c}(\yn)}\\
 &&~~~< \left\{\begin{array}{ll} C\displaystyle\sum_{\alpha\in\left(\scrptA^{*}-\{\alpha_c\}\right)} \left(p^{1+R/(1-\epsilon)}/n^2\right)^{\pna/2}/\sqrt{\pi} & \mbox{for } \nu =1,\\
C\displaystyle\sum_{\alpha\in\left(\scrptA-\{\alpha_c\}\right)} \left(2p^{R/(1-\epsilon)}/n^2\right)^{\pna/2} & \mbox{for } \nu =p.\\
\end{array}\right.\\
&&~~~ <  \left\{\begin{array}{ll} C\displaystyle\left\{p^{1/2+R/\{2(1-\epsilon)\}}/\left(n-p^{1/2+R/\{2(1-\epsilon)\}}\right)  \right\} & \mbox{for } \nu =1, \\
 C\displaystyle\left\{ \sqrt{2}p^{R/\{2(1-\epsilon)\}}/\left(n - \sqrt{2}p^{R/\{2(1-\epsilon)\}}\right) \right\}& \mbox{for } \nu =p.\\
\end{array}\right.
\end{eqnarray*}
One can choose $R>0$ and $\epsilon>0$ suitably and show that the above quantities go to 0 for any $0<b<1$.  \qed

\subsection{Proof of Theorem \ref{thm_2}}
Let $M_{\alpha_c}=M_{N}$. From (\ref{a.20}) and (A.4) we have
\begin{eqnarray}
 \sum_{\alpha\in\scrptA-\{\alpha_c\}} \frac{p(\Ma)}{p(\Mac)}\frac{m_{\alpha}(\yn)}{m_{\alpha_c}(\yn)}\geq \frac{1}{C} \sum_{\alpha}\int_{0}^{\infty} (1+g)^{-\pna/2} \pi(g) dg \label{a.40}
\end{eqnarray}
where $\pi(g)$ is given by (\ref{14}). Putting the prior we get the R.H.S. of the above expression as
\begin{eqnarray*}
 \frac{1}{C}\sum_{\alpha\in\scrptA-\{\alpha_c\}} \frac{\Gamma{(\gamma_0+\gamma_1)}\Gamma{(\gamma_1+\pna/2)}}{\Gamma{\gamma_1}\Gamma{(\gamma_0+\gamma_1+\pna/2)}} .
\end{eqnarray*}
Using the inequality of \cite{Gm_bound} stated above and the fact that $\gamma_1>\epsilon$ for some $\epsilon>0$ free of $n$, it can be shown that for some constant $C^{\prime}>0$, the above expression is bigger than
\begin{eqnarray*}
  C^{\prime}\sum_{\alpha\in\scrptA-\{\alpha_c\}} \left( \frac{\gamma_1}{\gamma_0+\gamma_1}\right)^{\pna/2}
  = C^{\prime} \left\{\left( 1+\sqrt{\frac{\gamma_1}{\gamma_0+\gamma_1}}\right)^{p}-1\right\}.
\end{eqnarray*}
Thus if $p=n^b$, the R.H.S. of (\ref{a.40}) does not go to $0$ if $\gamma_0/\gamma_1 = O(n^b)$. \qed

\subsection{Proof of Theorem \ref{thm_5}}
Our model selection criterion is to choose a model $\hat{\alpha}$ in the model space $\scrptA$, which maximizes $p(\Ma) m_{\alpha} (\yn)$ with respect to $\alpha$. Now, from Result \ref{res_1}, this is equivalent to maximizing
$$ p(\Ma) \Gamma{\left(\frac{\nu+\pna}{2}\right)} \{n\sysq(1-\rasq)\}^{-(n-1)/2} \left( \frac{n^2\nu}{2}\right)^{-\pna/2} (1+\varepsilon_n(\alpha)),$$
where $|\varepsilon_n(\alpha)|= p O_p(1) /(n \nu)$ uniformly in $\alpha$. We omit the other terms involved in the approximation of Result \ref{res_1}, since those are free of $\alpha$. Maximizing the above is equivalent to minimizing
\begin{equation}
\left[ p(\Ma) \Gamma{\left(\frac{\nu+\pna}{2}\right)} (1+\varepsilon_n(\alpha))\right]^{-2/(n-1)} \left( \frac{n^2\nu}{2}\right)^{\pna/(n-1)} n\sysq(1-\rasq)  \label{a.41}
\end{equation}
with respect to $\alpha$.
 From (\ref{20}) we have
$ n\sysq (1-\rasq) = C_n + 2 \sigma^2 D_n(\alpha) (1+\xi_n(\alpha)) $, where $C_n= \en^{\prime}\en$ and $\xi_n(\alpha)=\{2 \mun^{\prime}(I -\hatmat) \en - \en^{\prime} \hatmat \en \} /(2\sigma^2 D_n(\alpha))$. Hence, if $M_{\hat{\alpha}}$ is the model for which (\ref{a.41}) is minimized, then
\begin{equation*}
 \frac{D_n(\hat{\alpha})}{D_n(\alpha)} \leq \frac{C_n (b_{n}(\alpha) -1)}{2\sigma^2 D_n(\alpha)(1+\xi_n(\hat{\alpha}))} + \frac{b_n(\alpha)(1+\xi_n(\alpha)) }{(1+\xi_n(\hat{\alpha}))},
\end{equation*}
where 
\begin{align}
b_n(\alpha)= \left(\frac{P(M_{\hat{\alpha}})}{p(\Ma)}\right)^{2/(n-1)} \left(\frac{\Gamma{\{(n+p(\hat{\alpha}))/2\}}}{\Gamma{\{(n+\pna)/2\}}}\right)^{2/(n-1)} \left(\frac{n^2\nu}{2} \right)^{(p(\hat{\alpha}-\pna)/(n-1)}\nonumber\\ \left(\frac{1+\varepsilon_n(\hat{\alpha})}{1+\varepsilon_n(\alpha)} \right)^{2/(n-1)}. \label{a.42}
\end{align}
Therefore, if $\xi_n=\max_{\alpha} |\xi_n(\alpha)|$, we have
\begin{equation}
 1 \leq \frac{D_n(\hat{\alpha})}{\min_{\alpha} D_n(\alpha)} \leq \frac{C_n }{2n\sigma^2 (1-\xi_n)} \times \displaystyle\max_{\alpha} \frac{n(b_{n}(\alpha) -1)}{ D_n(\alpha)} + \frac{(1+\xi_n) }{(1-\xi_n)} \times \displaystyle\max_{\alpha} b_n(\alpha).
\end{equation}
The rest of the proof will follow from the following facts
\begin{eqnarray}
&& C_n/n \xrightarrow{p} \sigma^2, \label{a.43}\\
&& \xi_n \xrightarrow{p} 0,  \label{a.44} \\
&& \displaystyle\max_{\alpha} \frac{n(b_{n}(\alpha) -1)}{ D_n(\alpha)} \xrightarrow{p} 0, \label{a.45} \\
\mbox{and}\quad&& \displaystyle\max_{\alpha} b_n(\alpha) \xrightarrow{p} 1. \label{a.46}
\end{eqnarray}
The proof of (\ref{a.43}) is straight forward. To prove (\ref{a.44}) we note that
$$\xi_n \leq \frac{2 \max_{\alpha} \mun^{\prime} (I - \hatmat ) \en/n - \min_{\alpha}\en^{\prime}\hatmat\en/n}{2 \min_{\alpha}\sigma^2 D_n(\alpha)/n} \leq \frac{O_p(\sqrt{p/n})}{\delta/n^s},$$
from (iv) and (v) of Lemma \ref{lm_2} and assumption (A.3). Clearly if $s <(1-b)/2$, (\ref{a.44}) holds.

Next we prove (\ref{a.46}). We show that $ \log (\max_{\alpha}b_n(\alpha)) = \max_{\alpha} \log (b_n(\alpha)) \xrightarrow{p} 0$. From (\ref{a.42}) we have
\begin{eqnarray*}
  \max_{\alpha} \log b_n(\alpha)  & \leq & \frac{2}{n-1} \left\{ \log \left(\max_{\alpha}\frac{P(M_{\hat{\alpha}})}{p(\Ma)}\right) + \log \left(\max_{\alpha}\frac{\Gamma{\{(n+p(\hat{\alpha}))/2\}}}{\Gamma{\{(n+\pna)/2\}}}\right) \right.\\ 
  &&\left. +\log \left(\frac{1+\varepsilon_n(\hat{\alpha})}{1-\max_{\alpha}|\varepsilon_n(\alpha)|} \right) \right\} + \max_{\alpha} \frac{p(\hat{\alpha})-\pna}{n-1} \log \left(\frac{n^2\nu}{2} \right) \\
  &\leq& \frac{2}{n-1} \left\{ \log C + \frac{p}{2}\log \left(\frac{p+\nu}{2}\right) + \log \left(\frac{1+p O_p(1) / (n\nu)}{1-pO_p(1)/(n\nu) }\right) \right\}\\
  &&\hspace{2.5 in}+ \frac{p}{n-1} \log \left(\frac{n^2\nu}{2} \right),
\end{eqnarray*}
form (\ref{a.16}) and assumption (A.4).
It is now easy to show that when $p=O(n^b), ~0<b<1$, the above expression is $O_p\left( n^{-(1-b)} \log(n)\right)$ and converges to 0 with probability tending to 1. Hence (\ref{a.46}) holds.

Finally we prove (\ref{a.45}). By mean value theorem, for some $z>0$, $(e^z-1)=z e^{z^*} < z e^z$, where $z^{*} \in [0,z]$. Replacing $z$ by $\log b_n(\alpha)$ we get
$$ \max_{\alpha} (b_n(\alpha)-1) \leq \max_{\alpha} \log b_n(\alpha) \exp \{{\max_{\alpha} \log b_n(\alpha)}\}.$$
Thus by assumption (A.3) we have
\begin{align*} \displaystyle\max_{\alpha} \frac{n(b_{n}(\alpha) -1)}{ D_n(\alpha)} \leq \frac{\max_{\alpha}(b_{n}(\alpha) -1)}{\min_{\alpha} D_n(\alpha)/n} \leq n^s O_p\left( n^{-(1-b)} \log(n)\right) \\
\exp \left\{ O_p\left( n^{-(1-b)} \log(n)\right) \right\},\end{align*}
which is going to 0, with probability tending to 1. \qed

\vskip10pt
\noindent{\bf Acknowledgements.} I am thankful to my Ph.D. supervisor Professor Tapas Samanta for his guidance. I am also tankful to Professor Arijit Chaktabarti for his insightful comments. Their suggestions have helped to improve this paper substantially.
\bibliographystyle{imsart-nameyear}
\bibliography{ref.bib}
\end{document}